\documentclass[12pt]{article}

\usepackage{tikz}
\usetikzlibrary{automata}
\usetikzlibrary{decorations.pathreplacing}
\usepackage{forest}

\usepackage{amsmath,amssymb,amsthm}
\usepackage{cancel}
\usepackage{thmtools}

\usepackage{algorithm}
\usepackage[noend]{algpseudocode}

\usepackage{todonotes}
\usepackage{longtable}
\usepackage[f]{esvect}

\usepackage[colorlinks=true,allcolors=.]{hyperref}
\usepackage[nameinlink]{cleveref}

\usepackage{color}

\usepackage{authblk}

\usepackage{enumitem}
\usepackage{makecell}

\newtheorem{theorem}{Theorem}[section]
\newtheorem{lemma}[theorem]{Lemma}
\newtheorem{proposition}[theorem]{Proposition}
\newtheorem{corollary}[theorem]{Corollary}
\newtheorem{claim}[theorem]{Claim}
\newtheorem{observation}[theorem]{Observation}
\newtheorem{remark}[theorem]{Remark}
\newtheorem{conjecture}[theorem]{Conjecture}

\theoremstyle{definition}
\newtheorem{definition}[theorem]{Definition}
\newtheorem{example}[theorem]{Example}
\AtBeginEnvironment{example}{%
  \pushQED{\qed}%
}
\AtEndEnvironment{example}{\popQED}

\newcommand{\Ret}{\mathcal{R}}

\DeclareMathOperator{\RT}{\text{\normalfont RT}}

\def\RT{\mathrm{RT}}

\def\N{\mathbb N}

\def\R{\mathcal R}
\def\uu{\mathbf{u}}
\def\zz{\mathbf{v}}
\def\xx{\mathbf{x}}
\def\yy{\mathbf{y}}
\def\zz{\mathbf{z}}
\def\ww{\mathbf{w}}
\DeclareMathOperator{\lcp}{lcp}
\DeclareMathOperator{\lcs}{lcs}

\renewcommand{\tt}[1]{\mathtt{#1}}

\begin{document}

\title{An automated approach to structure theorems for minimally repetitive rich words}

\date{}

\author[1]{L\!'ubom\'{i}ra Dvo\v{r}\'{a}kov\'{a}}
\author[2]{Lucas Mol}
\author[1]{Edita Pelantov\'{a}}

\affil[1]{Department of Mathematics, Faculty of Nuclear Sciences and 
Physical Engineering, Czech Technical University in Prague\\
Trojanova 13, 120 00 Praha 2, Czech Republic\\
\url{lubomira.dvorakova,edita.pelantova@fjfi.cvut.cz}}

\affil[2]{Department of Mathematics and Statistics, Thompson Rivers University\\
805 TRU Way, Kamloops, BC V2C 0C8, Canada\\
\url{lmol@tru.ca}}

\maketitle

\begin{abstract}
Recent work on repetitions in rich words has led to the intensive study of the family of infinite rich words $\{\uu_d\,:\,d\geq 2\}$, where $\uu_d$ is the infinite fixed point of the morphism $\varphi_d:\Sigma_d^*\rightarrow \Sigma_d^*$ defined by $\varphi_d(i)=\tt{0}(i+1)$ for all $i\in\{\tt{0},\tt{1},\ldots,d-2\}$ and $\varphi_d(d-1)=\tt{0}(d-1)(d-1)$.  For $d=2$ and $d=3$, it is known that the repetition threshold for infinite $d$-ary rich words is equal to the asymptotic critical exponent $E^*(\uu_{2d-1})$ of~$\uu_{2d-1}$.  Moreover, the structure of $\uu_{2d-1}$ underlies the structure of all ``minimally repetitive'' $d$-ary rich words in the following sense: there is a real number $t_d>E^*(\uu_{2d-1})$ such that every infinite $t_d$-power-free $d$-ary rich word contains the language of $\pi(\uu_{2d-1})$ for one of only finitely many morphisms $\pi:\Sigma_{2d-1}^*\rightarrow \Sigma_d^*$.  We conjecture that a similar situation holds for all $d\geq 2$, and prove this conjecture for the first unknown case $d=4$.
In order to do this, we develop a new approach for proving structure theorems of the form described above.  Two new key concepts are the notion of ``compatibility'' of return words to palindromes in recurrent 
rich words, and the notion of ``$t$-persistent repetitions''.   These new concepts facilitate the automation of several aspects of the proofs of the structure theorems which were previously done laboriously by hand.  We first demonstrate the utility of our new approach by providing a much shorter proof of the previously known structure theorem for infinite $\tfrac{16}{7}$-power-free ternary rich words.  Then we prove a brand new structure theorem for quaternary rich words, which says that every infinite $2.12$-power-free quaternary rich word contains the language of one of seven distinct words of the form $\pi(\uu_7)$, up to permutation of the letters.  We prove that all seven of these words are rich and have critical exponent equal to $E^*(\uu_7)$, from which it follows that the repetition threshold for quaternary rich words is equal to $E^*(\uu_7)\doteq 2.1197196907$.

\end{abstract}

\section{Introduction}\label{Section:Introduction}

Since the completion of the proof of Dejean's theorem~\cite{CurrieRampersad2011,Rao11}, which describes the boundary between avoidable and unavoidable repetitions in infinite words over a fixed alphabet, there has been much research on determining the analogous ``repetition threshold'' for more restricted classes of infinite words.  This includes work on episturmian words~\cite{DvorakovaPelantova2024}, balanced words~\cite{DolceDP2023,DvorakovaPelantovaOpocenskaShur2022,RampersadShallitVandomme2019}, rich words~\cite{CuMoPe2024,CurrieMolRampersad2020}, and words with given factor complexity~\cite{Currie2025,OllingerShallit2025,ShallitShur2019}.  In this paper, we continue the study of repetitions in infinite rich words. 

A finite word $w$ of length $n$ is called {\em rich} if it contains exactly $n+1$ distinct palindromes as factors, including the empty word.  Such words are ``rich'' in palindromes, as every word of length $n$ contains \emph{at most} $n+1$ distinct palindromes as factors.
An infinite word is {\em rich} if all of its (finite) factors are rich.
Since their introduction about two decades ago~\cite{DroubayJustinPirillo2001}, much work has been done on rich words from many different perspectives.
\begin{itemize}
\item Examples of infinite rich words include episturmian words~\cite{DroubayJustinPirillo2001}, complementary-symmetric Rote words~\cite{BBLV2011}, and words coding symmetric $d$-interval exchange~\cite{BaMaPe2007}.
\item There are characterizations of infinite rich words in terms of complete return words~\cite{GlenJustinWidmerZamboni2009}, a relation between factor and palindromic complexity~\cite{BuLuGlZa2009}, and bi-extensions of bispecial factors~\cite{BalkovaPelantovaStarosta2010}.
\item Upper and lower bounds on the number of rich words of length $n$ over an alphabet of size $d$ have been established~\cite{Guo2016,Rukavicka2017}.  An asymptotically accurate count is not yet known.
\item Morphisms preserving richness have been studied by several authors~\cite{BPS2011,DolcePelantova2022,GlenJustinWidmerZamboni2009}. A~complete characterization of such morphisms is not yet known.
\end{itemize}

Let us recall the essential notions related to repetitions. For a positive integer $n$, the $n$-th power of a nonempty word $u$, denoted $u^n$, is a repetition of $n$ copies of $u$.
In other words, $u^n$ denotes the prefix of length $|u|n$ of the infinite periodic word $uuu\cdots = u^\omega$, where $|u|$ is the length of $u$. 
More generally, we can define rational powers: the prefix of length $k$ of $u^\omega$ can be written as $u^e$, where $e=\frac{k}{|u|}$. 
For a real number $t>1$, we say that an infinite word $\uu$ is $t$-{\em power-free} if for every nonempty factor of $\uu$ of the form $u^e$, we have $e<t$.
The {\em critical exponent} of an infinite word $\uu$, denoted $E(\uu)$, is defined as
$$E(\uu) =\sup\{e \in \mathbb{Q}: \  u ^e \  \text{is a factor of   } \uu  \  \text{for a nonempty word} \  u\}\,.
$$
The {\em asymptotic critical exponent} of $\uu$, denoted by $E^*(\uu)$, was introduced by Cassaigne~\cite{Ca2008}, who called it the \emph{asymptotic index} of $\uu$. It is defined
to be $+\infty$,  if $E(\uu) = +\infty$, and otherwise by
 $$E^*(\uu) =\limsup_{n\to \infty}\{e \in \mathbb{Q}: \  u ^e \  \text{is a factor of  } \uu  \  \text{for some word }  u \ \text{of length} \  n  \}\,.$$
It is readily seen that $E^*(\uu)\leq E(\uu)$.

\medskip

We are interested in the boundary between avoidable and unavoidable repetitions in infinite rich words over $\Sigma_d=\{\tt{0},\tt{1},\ldots,d-1\}$.  This is the {\em repetition threshold} for the class $\mathcal{C}_d$ of infinite rich words over $\Sigma_d$, which we denote by $\RT(\mathcal{C}_d)$, and define formally as
$$\RT(\mathcal{C}_d) =\inf\{E(\uu): \uu \text{\ is an infinite rich word over $\Sigma_d$}\}\,.$$ 
The analogous \emph{asymptotic repetition threshold} for infinite rich words over $\Sigma_d$ will also be relevant to our discussion.  It is denoted by $\RT^*(\mathcal{C}_d)$, and defined by 
$$\RT^*(\mathcal{C}_d) =\inf\{E^*(\uu): \uu \text{\ is an infinite rich word over $\Sigma_d$}\}\,.$$
Since $E^*(\uu)\leq E(\uu)$ for every infinite word $\uu$, we have $\RT^*(\mathcal{C}_d)\leq \RT(\mathcal{C}_d)$.  

Further, since every recurrent rich word contains infinitely many distinct overlaps as factors~\cite{PStarosta2013}, it follows that $\RT(\mathcal{C}_d)\geq \RT^*(\mathcal{C}_d)\geq 2$.

One particular family of infinite words appears to be at the core of the problem of determining $\RT(\mathcal{C}_d)$.  For each $d\geq 2$, we let $\uu_d$ denote the infinite fixed point of the morphism $\varphi_d:\Sigma_d^*\rightarrow \Sigma_d^*$ defined by
\begin{equation}\label{def:varphi} 
\begin{split}
   \varphi_d(i)&=\tt{0}(i+1)\quad  \text{for $i\in\{\tt{0},\ldots,d-2\}$}\,, \text{ and } \\
   \varphi_d(d-1)&=\tt{0}(d-1)(d-1)\,.
\end{split}
\end{equation}
For example, with $d=5$, we have
\begin{align*}
    \varphi_5(\tt{0})&=\tt{01}\\
    \varphi_5(\tt{1})&=\tt{02}\\
    \varphi_5(\tt{2})&=\tt{03}\\
    \varphi_5(\tt{3})&=\tt{04}\\
    \varphi_5(\tt{4})&=\tt{044}
\end{align*}
and
\[
\uu_5=\tt{0102010301020104010201030102010440102010\cdots}\,.
\]
Dvo\v{r}\'{a}kov\'{a}, Klouda, and Pelantov\'{a}~\cite{DvorakovaKloudaPelantova2025} have shown that for all $d\geq 2$, the word $\uu_d$ is rich and satisfies
\[
E^*(\uu_d)=1+\tfrac{1}{3-\Lambda_d},
\]
where $\Lambda_d$ is the spectral radius of the incidence matrix of $\varphi_d$.
This knowledge of $\uu_d$ has been leveraged to establish the following general results~\cite{DvorakovaKloudaPelantova2025,DvPe_LimRT2025}.
\begin{itemize}
    \item For each $d\geq 2$, there is a morphism $\rho_d:\Sigma_d^*\rightarrow \Sigma_2^*$ such that $\rho_d(\uu_d)$ is rich and $E^*(\rho_d(\uu_d))=E^*(\uu_d)$.  It follows that $\RT^*(\mathcal{C}_d)=2$ for all $d\geq 2$, as $\displaystyle\lim_{d\to \infty }\Lambda_d = 2$.
    \item For each $d\geq 1$, there is a morphism $\sigma_d:\Sigma_{2d+1}^*\rightarrow \Sigma_{d+2^d}^*$ such that $\sigma_d(\uu_{2d+1})$ is rich and satisfies $E(\sigma_d(\uu_{2d+1}))=E^*(\uu_{2d+1})$.  It follows that $\displaystyle\lim_{d\rightarrow \infty} \RT(\mathcal{C}_d)=2$. 
\end{itemize}

We are interested in the exact values of $\RT(\mathcal{C}_d)$, which are currently known only for $d=2$ and $d=3$.  Currie, Mol, and Rampersad~\cite{CurrieMolRampersad2020} showed that
\[
\RT(\mathcal{C}_2)=E^*(\uu_3)\doteq 2.70710678,
\]
and Currie, Mol, and Peltom\"aki~\cite{CuMoPe2024} showed that
\[
\RT(\mathcal{C}_3)=E^*(\uu_5)\doteq 2.25876324.
\]
These two results lead naturally to the following conjecture.

\begin{conjecture}[\cite{DvPe_LimRT2025}]\label{con:opet}  
    For every integer $d\geq 2$, we have $\RT(\mathcal{C}_{d}) = E^*(\uu_{2d-1})$.
\end{conjecture}

Let us describe the apparent connection between $\RT(\mathcal{C}_d)$ and the word $\uu_{2d-1}$ in more detail by explaining how the values of $\RT(\mathcal{C}_2)$ and $\RT(\mathcal{C}_3)$ were established.  The key results are the two structure theorems stated below.\footnote{We have modified the original statements of the structure theorems for ease of reading and to align with the more complete picture of the general situation that we now have, but we stress that \autoref{Theorem:BinaryStructure} and \autoref{Theorem:TernaryStructure} follow directly from the results referenced at the beginning of their statements.}  
Roughly speaking, for $d=2$ and $d=3$, respectively, they say that every infinite rich word over $d$ letters with critical exponent sufficiently close to $E^*(\uu_{2d-1})$ must ``look like'' the image of $\uu_{2d-1}$ under some morphism $\pi$ from a~finite set of morphisms $\mathcal{P}_d$.  Moreover, for each $\pi\in\mathcal{P}_d$, the word $\pi(\uu_{2d-1})$ is rich and has critical exponent equal to $E^*(\uu_{2d-1})$.  It follows directly from these statements that $\RT(\mathcal{C}_d)=E^*(\uu_{2d-1})$ for $d=2$ and $d=3$.  Since the morphisms $\pi\in\mathcal{P}_d$ map words over the larger alphabet $\Sigma_{2d-1}$ to the smaller alphabet $\Sigma_d$, we call them \emph{projections}.

\begin{theorem}[c.f.~{\cite[Theorem 1]{CurrieMolRampersad2020}}]\label{Theorem:BinaryStructure}
    If $\zz$ is an infinite $\frac{14}{5}$-power-free binary rich word, then $\zz$ contains every factor of either $\pi_1(\uu_3)$ or $\pi_2(\uu_3)$, up to permutation of the letters, where $\pi_1,\pi_2:\Sigma_3^*\rightarrow \Sigma_2^*$ are the projections defined below.\footnote{In the work of Currie, Mol, and Rampersad~\cite{CurrieMolRampersad2020}, the morphism $\pi_1$ defined here is denoted by $f$, the morphism $\pi_2$ defined here is denoted by $f\circ {g}$, and the morphism $\varphi_3$ is denoted by $h$, so that $\uu_3=h^\omega(\tt{0})$.}
    \begin{center}
    \begin{tabular}{lll}
$\pi_1: \begin{array}{l}
    \tt{0} \mapsto \tt{0}\\
    \tt{1} \mapsto \tt{01}\\
    \tt{2} \mapsto \tt{011}
\end{array}$ & \quad & 
$\pi_2: \begin{array}{l}
    \tt{0} \mapsto \tt{00101}\\
    \tt{1} \mapsto \tt{00101101}\\
    \tt{2} \mapsto \tt{0010110101101}
\end{array}$
\end{tabular}
\end{center}
Moreover, the words $\pi_1(\uu_3)$ and $\pi_2(\uu_3)$ are rich and have critical exponent equal to $E^*(\uu_3)<\frac{14}{5}$.
\end{theorem}

\begin{theorem}[c.f.~{\cite[Theorem 1.2]{CuMoPe2024}} and~{\cite[Proposition 10]{DvPe2025}}]\label{Theorem:TernaryStructure}
    If $\zz$ is an infinite $\tfrac{16}{7}$-power-free ternary rich word, then $\zz$ contains every factor of $\pi(\uu_5)$, up to permutation of the letters, where $\pi:\Sigma_5^*\rightarrow \Sigma_3^*$ is the projection defined below.\footnote{This structure theorem was proven by Currie, Mol, and Peltom\"{a}ki~\cite{CuMoPe2024}, who used an alternate construction of the word $\pi(\uu_5)$.  The morphic construction $\pi(\uu_5)$ was found by Dvo\v{r}\'{a}kov\'{a} and Pelantov\'{a}~\cite{DvPe2025}.}
    \[
    \pi: \begin{array}{l}
    \tt{0} \mapsto A\\
    \tt{1} \mapsto AD\\
    \tt{2} \mapsto AC\\
    \tt{3} \mapsto ACC\\
    \tt{4} \mapsto ACCBCC,
\end{array}
    \]
    where $A=\tt{00101101}$, $B=\tt{001}$, $C=\tt{00202202}$, and $D=\tt{002}$.  Moreover, the word $\pi(\uu_5)$ is rich and has critical exponent equal to $E^*(\uu_5)<\tfrac{16}{7}$.
\end{theorem}

We also note that the words $\pi_1(\uu_3)$, $\pi_2(\uu_3)$, and $\pi(\uu_5)$ are not only rich in palindromes, but exhibit further symmetry.  In fact, each of these words is \emph{$G$-rich} (see~\cite{PStarosta2013} for the definition) for a group $G$ generated by two involutive antimorphisms~\cite{DvPe2025}.  In each case, one of these antimorphisms is the reversal antimorphism, and the other swaps certain pairs of letters.  In particular, this means that the words $\pi_1(\uu_3)$, $\pi_2(\uu_3)$, and $\pi(\uu_5)$ exhibit strong symmetry under both of these antimorphisms.

All of this leads us to the following conjecture, which generalizes the structure theorems stated above for $d=2$ and $d=3$, and from which \autoref{con:opet} follows easily.

\begin{conjecture}\label{con:structure}
    For every integer $d\geq 2$, there is some real number $t_d>E^*(\uu_{2d-1})$ and some finite set $\mathcal{P}_d$ of projections from $\Sigma_{2d-1}^*$ to $\Sigma_d^*$ such that if $\zz$ is an infinite $t_d$-power-free rich word over $d$ letters, then 
    $\zz$ contains every factor of $\pi(\uu_{2d-1})$ for some $\pi\in\mathcal{P}_d$.  Moreover, for each projection $\pi\in\mathcal{P}_d$, the word $\pi(\uu_{2d-1})$ is rich and has critical exponent $E^*(\uu_{2d-1})$, and is $G$-rich for some group $G$ generated by two involutive antimorphisms.
\end{conjecture}

In this paper, we develop an automated approach to proving the structural part of \autoref{con:structure} for a fixed $d$.  We demonstrate the utility of this new approach by giving a shorter proof of the structural part of \autoref{Theorem:TernaryStructure}.  Then we provide a complete proof of \autoref{con:structure} in the case $d=4$ by establishing the theorem below.

\begin{theorem}\label{Theorem:QuaternaryStructure}
    If $\zz$ is an infinite $\tfrac{53}{25}$-power-free quaternary rich word, then up to permutation of the letters, $\zz$ contains every factor of $\pi_i(\uu_7)$ for some $i\in\{1,2,\ldots,7\}$, where the $\pi_i$'s are defined in \autoref{Table:PiDefinition}.  Moreover, for each $i\in\{1,2,\ldots,7\}$, the word $\pi_i(\uu_7)$ is rich and has critical exponent $E^*(\uu_7)<\tfrac{53}{25}$, and is $G$-rich for some group $G$ generated by two involutive antimorphisms.
\end{theorem}

The following corollary is immediate.
\begin{corollary}
    $\RT(\mathcal{C}_4)= E^*(\uu_7)\doteq 2.1197196907$.
\end{corollary}

Let us briefly summarize our overall approach on the structural side.  Once we have candidates for $t_d$ and $\mathcal{P}_d$, we wish to show that if $\zz$ is an infinite $t_d$-power-free rich word over $d$ letters, then $z$ contains every factor of $\pi(\uu_{2d-1})$ for some $\pi\in\mathcal{P}_d$.  Recall that $\uu_{2d-1}=\varphi^\omega(\tt{0})$, where we write $\varphi=\varphi_{2d-1}$ for ease of notation.
Roughly speaking, the idea is to ``peel back the layers'' of structure.  For each infinite $t_d$-power-free quaternary rich word $\zz$, we show that 
\begin{itemize}
    \item $\zz$ has a suffix of the form $\pi(\yy)$ for some $\pi\in\mathcal{P}_d$,
    \item $\yy$ has a suffix of the form $\varphi(\xx)$ (so that $\zz$ has a suffix of the form $\pi(\varphi(\xx))$),
    \item $\xx$ has a suffix of the form $\varphi(\ww)$ (so that $\zz$ has a suffix of the form $\pi(\varphi^2(\ww))$),
    \item and so on, so that for all $n\geq 0$, the word $\zz$ has a suffix of the form $\pi(\varphi^n(\mathbf{v}))$.
\end{itemize} 
Establishing each ``layer'' of structure indicated above involves determining all possible return words to an appropriate palindrome in the given word, and observing that they are exactly the images of the letters under the desired morphism.  We refer to the first layer as the ``outer layer'' (which will serve as the base case of an inductive argument), and to the remaining layers as ``inner layers'' (which will be covered by the inductive step).\footnote{For technical reasons, we actually replace each $\pi\in\mathcal{P}_d$ in the argument above with $\hat{\pi}=\pi\circ \varphi^k$ for some appropriate positive integer $k$, depending on $d$.  Essentially, we include a small number of iterations of $\varphi$ in the ``outer layer''.  This is necessary to get the inductive argument started.}

We develop two new key concepts which facilitate our automated approach to establishing these layers of structure.
\begin{itemize}
\item The first is the notion of ``compatibility''.  We prove what we call the \nameref{Lemma:Compatibility}, which says that any two return words to the same palindrome in a~recurrent rich word must satisfy a simple condition.  This condition is easy to check by computer, and if we know enough about a recurrent rich word $\ww$, then we can often use it to find all possible return words to a given palindrome in $\ww$ purely by computer.  (Remember that this is exactly how we establish each layer of structure.)

\item The second is the notion of ``$t$-persistent repetitions'' (with respect to the morphisms $\varphi$ and $\pi$).  Roughly speaking, these are factors in a word $\ww$ that lead to factors of exponent at least $t$ in $\pi(\varphi^n(\ww))$ for all $n\geq 0$.  Evidently, if $\pi(\varphi^k(\ww))$ is $t_d$-power-free for some $k\geq 0$, then $\ww$ must not contain any $t_d$-persistent repetitions.  We develop an efficient method for verifying that certain repetitions are $t$-persistent by computer, which we call the \nameref{Proposition:Persistence}.  This is the key to automating the proof of the inner layers of structure.
\end{itemize}

\begin{table}
\begin{longtable}{lll}
$\pi_1: \begin{array}{l}
    \tt{0} \mapsto \tt{0}\\
    \tt{1} \mapsto \tt{01}\\
    \tt{2} \mapsto \tt{011}\\
    \tt{3} \mapsto \tt{02}\\
    \tt{4} \mapsto \tt{022}\\
    \tt{5} \mapsto \tt{022322}\\
    \tt{6} \mapsto \tt{0223223322322}, 
\end{array}$ & \quad & 
$\pi_2: \begin{array}{l}
    \tt{0} \mapsto \tt{01}\\
    \tt{1} \mapsto \tt{0}\\
    \tt{2} \mapsto \tt{011}\\
    \tt{3} \mapsto \tt{02}\\
    \tt{4} \mapsto \tt{0232}\\
    \tt{5} \mapsto \tt{0232232}\\
    \tt{6} \mapsto \tt{023223233232232},
\end{array}$\\ & &\\
$\pi_3: \begin{array}{l}
    \tt{0} \mapsto \tt{01}\\
    \tt{1} \mapsto \tt{011}\\
    \tt{2} \mapsto \tt{0}\\
    \tt{3} \mapsto \tt{02}\\
    \tt{4} \mapsto \tt{0232}\\
    \tt{5} \mapsto \tt{023233232}\\
    \tt{6} \mapsto \tt{02323323223233232},
\end{array}$ &\qquad & 
$\pi_4: 
\begin{array}{l}
    \tt{0} \mapsto \tt{0110}\\
    \tt{1} \mapsto \tt{01102}\\
    \tt{2} \mapsto \tt{0110232}\\
    \tt{3} \mapsto \tt{01102332}\\
    \tt{4} \mapsto \tt{0110}A\\
    \tt{5} \mapsto \tt{0110}A\tt{0}A\\
    \tt{6} \mapsto \tt{0110}A\tt{0}A\tt{010}A\tt{0}A,\\
    \text{\quad  \ where $A=\tt{23322332}$,}
\end{array}$\\ && \\
 $\pi_5: \begin{array}{l}
    \tt{0} \mapsto \tt{01011}\\
    \tt{1} \mapsto \tt{01022}\\
    \tt{2} \mapsto \tt{0102}\\
    \tt{3} \mapsto \tt{010232}\\
    \tt{4} \mapsto \tt{010}B\\
    \tt{5} \mapsto \tt{010}B\tt{00}B\\
    \tt{6} \mapsto \tt{010}B\tt{00}B\tt{0}B\tt{00}B,\\
    \text{\quad  \ where $B=\tt{23233232}$,}\end{array}$
 &\quad & 
 $\pi_6: \begin{array}{l}
    \tt{0} \mapsto \tt{00101101002}\\
    \tt{1} \mapsto \tt{001011010022}\\
    \tt{2} \mapsto \tt{001011010022322}\\
    \tt{3} \mapsto \tt{0010110100}C\\
    \tt{4} \mapsto \tt{0010110100}C\tt{0}C\\
    \tt{5} \mapsto \tt{0010110100}C\tt{0}C\tt{00}C\tt{0}C\\
    \tt{6} \mapsto \tt{0010110100}C\tt{0}C\tt{00}C\tt{0}C\tt{00100}C\tt{0}C\tt{00}C\tt{0}C,\\
    \text{\quad \ where $C=\tt{2232332322}$,}\end{array}$\\ \\
    
\multicolumn{3}{l}{$\pi_7: \begin{array}{l}
    \tt{0} \mapsto \tt{0101100110102}\\
    \tt{1} \mapsto \tt{010110011010232}\\
    \tt{2} \mapsto \tt{01011001101023233232}\\
    \tt{3} \mapsto \tt{010110011010}D\\
    \tt{4} \mapsto \tt{010110011010}D\tt{0}D\\
    \tt{5} \mapsto \tt{010110011010}D\tt{0}D\tt{010}D\tt{0}D\\
    \tt{6} \mapsto \tt{010110011010}D\tt{0}D\tt{010}D\tt{0}D\tt{01011010}D\tt{0}D\tt{010}D\tt{0}D,\\
    \text{\quad \  where $D=\tt{232332233232}$.}
\end{array}$ }
\end{longtable}
\setcounter{table}{0}
\caption{The seven projections $\pi_1,\pi_2,\ldots,\pi_7$ of \autoref{Theorem:QuaternaryStructure}.}
\label{Table:PiDefinition}
\end{table}

The layout of the remainder of the article is as follows.  In \autoref{Section:Preliminaries}, we give the necessary background and preliminaries.  
In \autoref{Section:UsefulLemma}, we prove the \nameref{Lemma:Compatibility} for return words to palindromes in recurrent rich words.
In \autoref{Section:PersistentRepetitions}, we introduce $t$-persistent repetitions and prove the \nameref{Proposition:Persistence}.
In \autoref{Section:TernaryStructure}, we demonstrate the utility of our automated approach by giving a much shorter proof of the known structure theorem for $\tfrac{16}{7}$-power-free ternary rich words (i.e., the structural part of \autoref{Theorem:TernaryStructure}).
In \autoref{Section:QuaternaryStructure}, we prove our new structure theorem for $\tfrac{53}{25}$-power-free quaternary rich words (i.e., the structural part of \autoref{Theorem:QuaternaryStructure}).  In \autoref{Section:Pi1}, we complete the proof of \autoref{Theorem:QuaternaryStructure} by establishing the desired properties of $\pi_i(\uu_7)$ for all $i\in\{1,2,\ldots,7\}$.  We give a detailed proof for the case $i=1$ only, as the proofs of the other cases are analogous. 
Note that the structural part of \autoref{Theorem:TernaryStructure} together with the proof that $\pi_1(\uu_7)$ is rich and satisfies  $E(\pi_1(\uu_7)) = E^*(\uu_7)$ already confirms \autoref{con:opet} for $d=4$, as $E(\pi_i(\uu_7))\geq E^*(\pi_i(\uu_7))= E^*(\uu_7)$ for every $i=1,2,\ldots,7$.
Finally, in \autoref{Section:OpenQuestions}, we close with some open questions.

The proofs of the structural parts of \autoref{Theorem:TernaryStructure} and \autoref{Theorem:QuaternaryStructure}, and several examples throughout the paper, rely on extensive computations.  All of the code that we used to verify our results can be found at the following link.
\begin{center}
\url{https://github.com/lgmol/Structure-Theorems-for-Rich-Words/}
\end{center}
The printed output and run times of the files that we ran in order to verify our results can also be found there.

\section{Preliminaries}\label{Section:Preliminaries}

Throughout, for every integer $d\geq 2$, we let $\Sigma_d$ denote the \emph{alphabet} $\{\tt{0},\tt{1}, \ldots, d-1\}$.  The elements of $\Sigma_d$ are called \textit{letters}.  The symbol $\Sigma_d^*$ denotes the set of \textit{finite words} over $\Sigma_d$, i.e., the set of finite strings $w$ of the form $w=w_0w_1\cdots w_{n-1}$, where $n\in\N$ and $w_i \in \Sigma_d$ for all $i\in\{0,1,\ldots, n-1\}$.  For a word $w\in\Sigma_d^*$, the \emph{length} of $w$ is the number of letters that make up $w$, and is denoted by $|w|$.  There is a unique word of length $0$, which we denote by $\varepsilon$, and call the \emph{empty word}.  The set $\Sigma_d^*$ equipped with concatenation as the operation forms a monoid with the empty word $\varepsilon$ as its identity element.  For each $i\in\Sigma_d$, the number of occurrences of the letter $i$ in $w$ is denoted by $|w|_i$. The \textit{Parikh vector} of $w$, denoted by $\vv{w}$, is the vector whose $i$-th component is $|w|_{i}$, i.e., $\vv{w}=\left(|w|_\tt{0},|w|_\tt{1},\ldots,|w|_{d-1}\right)^T$. 
If we can write $w=pvs$ for (possibly empty) words $p,v,s\in\Sigma_d^*$, then the word $p$ is called a \textit{prefix} of $w$, the word $s$ is called a \textit{suffix} of $w$, and the word $v$ is called a \textit{factor} of $w$.  Two words $v,w$ are said to be \textit{conjugate}  if there exist words $x,y$ such that $w=xy$ and $v=yx$.

An~\textit{infinite word} $\ww$ over $\Sigma_d$ is an infinite string of the form $\ww = w_0 w_1 w_2 \cdots$, where $w_j \in \Sigma_d$ for all $j \in \N$.  We let $\Sigma_d^{\omega}$ denote the set of infinite words over $\Sigma_d$, and use boldface letters to denote infinite words.  A finite word $v$ over $\Sigma_d$ is called a~\textit{factor} of the infinite word $\ww = w_0 w_1 w_2 \cdots$ if there exists some $j \in \mathbb N$  such that $v= w_j w_{j+1}\cdots w_{j+|v|-1}$. In this case, the integer $j$ is called an \textit{occurrence} of the factor $v$ in the infinite word $\ww$. If there is an occurrence of $v$ at index $0$, then $v$ is called a \textit{prefix} of $\ww$. A \textit{suffix} of $\ww$ is an infinite word of the form $\mathbf{v}=w_jw_{j+1}w_{j+2}\cdots$ for some $j\in\N$.  

The set of all (finite) factors of an infinite word $\ww$ is called the \textit{language} of $\ww$, and is denoted by $\mathcal{L}(\ww)$. The set of factors of $\ww$ of length $n$ is denoted by $\mathcal{L}_n(\ww)$.
A~factor $w$ of an infinite word $\ww$ is called \textit{left special} if $iw, jw \in \mathcal{L}(\uu)$ for at least two distinct letters $i, j\in \Sigma_d$. A~\textit{right special} factor is defined analogously, and a factor is called \textit{bispecial} if it is both left and right special.  We define $\mathrm{Lext}(w)=\{iw\in\mathcal{L}(\ww)\,\colon\ i\in\Sigma_d\}$, and call the words in $\mathrm{Lext}(w)$ the \emph{left extensions} of $w$ in~$\ww$.  Similarly, we define $\mathrm{Rext}(w)=\{wj\in\mathcal{L}(\ww)\,\colon\, j\in\Sigma_d\}$, and call the words in $\mathrm{Rext}(w)$ the \emph{right extensions} of $w$ in~$\ww$.  Finally, we define $\mathrm{Bext}(w)=\{iwj\in\mathcal{L}(\ww)\,\colon\, i,j\in\Sigma_d\}$, and call the words in $\mathrm{Bext}(w)$ the \emph{bi-extensions} of $w$ in~$\ww$.
The \textit{bilateral order}~\footnote{The bilateral order was introduced by Cassaigne~\cite{Cassaigne1997} as a tool for the computation of factor complexity.} of $w$, denoted by $b(w)$, is defined by
$$\mathrm{b}(w)=\#\mathrm{Bext}(w)-\#\mathrm{Lext}(w)-\#\mathrm{Rext}(w)+1\,.$$


An infinite word $\ww$ is said to be \textit{recurrent} if for each $w\in\mathcal{L}(\ww)$, there are at least two occurrences of $w$ in $\ww$.  A recurrent word $\ww$ is said to be \textit{uniformly recurrent} if for each $w\in\mathcal{L}(\ww)$, the distances between consecutive occurrences of $w$ in $\ww$ are bounded by some constant (depending on $w$).
An infinite word $\ww$ is \textit{eventually periodic} if there exist words $x \in \Sigma_d^*$ and $v \in \Sigma_d^* \setminus \{\varepsilon\}$ such that $\ww$ can be written as $\ww = xvvv \cdots = xv^\omega$. If $\ww$ is not eventually periodic, then $\ww$ is called \textit{aperiodic}.

Consider a factor $v$ of a recurrent word $\ww = w_0 w_1 w_2 \cdots$. Let $j < \ell$ be two consecutive occurrences of $v$ in $\ww$. Then the word $u=w_j w_{j+1} \cdots w_{\ell-1}$ is called a \textit{return word} to $v$ in $\ww$, and the word $uv$ is called a \textit{complete return word} to $v$ in $\ww$.
The set of return words to $v$ in $\ww$ is denoted by ${\mathcal R}_{\ww}(v)$. Note that the set $\R_\ww(v)$ is finite for every $v \in \mathcal{L}(\ww)$ if and only if $\ww$ is uniformly recurrent.

A \textit{morphism} is a map $\psi: \Sigma_d^* \to \Sigma_k^*$ such that $\psi(uv) = \psi(u)\psi(v)$ for all words $u, v \in \Sigma_d^*$.
The morphism $\psi$ can be extended naturally to every infinite word $\ww=w_0 w_1 w_2\cdots$ over $\Sigma_d$ by setting
$\psi(\ww) = \psi(w_0) \psi(w_1) \psi(w_2) \cdots\,$.  A \textit{fixed point} of a morphism $\psi:  \Sigma_d^* \to \Sigma_d^*$ is an infinite word $\ww$ such that $\psi(\ww) = \ww$.
To each morphism $\psi:  \Sigma_d^* \to \Sigma_k^*$, we associate the $k\times d$ \textit{incidence matrix} $M_\psi$, defined by $(M_\psi)_{ij}=|\psi(j)|_{i}$ for all $i \in \Sigma_k$ and $j \in \Sigma_d$. 
We say that $\psi$ is \textit{primitive} if the matrix $M_\psi$ is primitive, i.e., if there exists $\ell\in \mathbb N$ such that $M_\psi^\ell$ has all positive entries.  

An \textit{antimorphism} is a map $\Theta: \Sigma_d^* \to \Sigma_d^*$ such that $\Theta(uv) = \Theta(v)\Theta(u)$ for all words $u, v \in \Sigma_d^*$.
The \textit{reversal} mapping $R:\Sigma_d^* \to \Sigma_d^*$ is the antimorphism satisfying $R(i)=i$ for each letter $i \in \Sigma_d$.  For each word $w\in\Sigma_d^*$, the word $R(w)$ is called the \emph{mirror image} or \emph{reversal} of $w$.
The language of $\ww$ is said to be \textit{closed under the antimorphism $\Theta$} if for each factor $w\in\mathcal{L}(\ww)$, the word $\Theta(w)$ is also a factor of $\ww$.
In particular, the language of $\ww$ is said to be \textit{closed under reversal} if for each factor $w\in\mathcal{L}(\ww)$, its mirror image $R(w)$ is also a factor of $\ww$. 
 
A word $p$ is called a \textit{palindrome} (or is said to be \emph{palindromic}) if it is equal to its own mirror image, i.e., if $p=R(p)$.
If a uniformly recurrent word $\ww$ contains infinitely many palindromic factors, then its language ${\mathcal L}(\ww)$ is closed under reversal. 
If $p$ is a~palindromic factor of an infinite word $\ww$ over $\Sigma_d$, then we let $\mathrm{Pext}(p)=\{ipi \in {\mathcal L}(\ww)\, :\, i \in {\Sigma_d}\}$, and call the words in $\mathrm{Pext}(p)$ the \emph{palindromic extensions} of $p$ in $\ww$.

Every finite word $u$ contains at most $|u|+1$ distinct palindromic factors (including the empty word). So if $u$ contains exactly $|u|+1$ distinct palindromic factors, then $u$ is said to be {\em rich} (in palindromes). An infinite word $\ww$ is said to be \emph{rich} if all of its (finite) factors are rich. 
We use the following equivalent characterization of richness~\cite{GlenJustinWidmerZamboni2009}. 

\begin{theorem}\label{Theorem:CompleteReturnsToPalindromes}
An infinite word $\ww$ is rich if and only if every complete return word to a palindromic factor of $\ww$ is a~palindrome.
\end{theorem}

Furstenberg~\cite{Furstenberg1981} proved that for every infinite word $\ww$ over an alphabet $\Sigma_d$, there exists a~uniformly recurrent  word $\uu$ such that ${\mathcal L}(\uu) \subset {\mathcal L}(\ww)$, which implies that 
$E(\uu)\leq E(\ww)$.  Recall that the language of finite rich words is \emph{factorial}, i.e., every factor of every rich word is rich.   So if $\ww$ is rich and ${\mathcal L}(\uu) \subset {\mathcal L}(\ww)$, then $\uu$ is rich as well. 
Consequently, when searching for the repetition threshold of the class ${\mathcal C}_d$ of infinite $d$-ary rich words, it suffices to consider uniformly recurrent words.

\begin{definition}[\cite{BPS2011}] \label{eq:Pret}
A morphism $\eta: \Sigma_d^* \to \Sigma_{d'}^*$ is said to be of class $P_{\mathrm{ret}}$ if $\eta$ is injective on the alphabet $\Sigma_d$ and there exists a~palindromic marker $p \in \Sigma_{d'}^*$ such that 
for all $i \in \Sigma_d$, the word $\eta(i)p$ is a~palindrome and it is a~complete return word to $p$. 
\end{definition}
Morphisms of class $P_{\mathrm{ret}}$ have been studied extensively in the context of rich words~\cite{BPS2011,DolcePelantova2022}, and play an important role in this paper.
\begin{remark}
For all $d\geq 2$, the morphism $\varphi_d$ belongs to class $P_{\mathrm{ret}}$ with marker $\tt{0}$.  Similarly, the projections $\pi_1,\pi_2,\ldots,\pi_7$ defined in \autoref{Table:PiDefinition}, the projections $\pi_1, \pi_2$ from \autoref{Theorem:BinaryStructure}, and the projection $\pi$ from \autoref{Theorem:TernaryStructure} all belong to class $P_{\mathrm{ret}}$.
\end{remark}

Let us recall some important results from~\cite{BPS2011} related to morphisms of class $P_{\mathrm{ret}}$.
\begin{itemize}
\item If an infinite rich word $\ww$ is the image of an infinite word $\uu$ under some morphism~$\eta$ of class $P_{\mathrm{ret}}$, then $\uu$ is rich, too. 
\item The image of an infinite rich word $\uu$ by a morphism $\eta \in P_{\mathrm{ret}}$ may not be rich. 
\item If $p$ is a palindromic factor of a~uniformly recurrent rich word~${\ww}$ and $r_0, r_1, \ldots, r_{d-1}$ is the list of return words to $p$ in $\ww$, then  a~suffix ${\ww'}$ of ${\ww}$ is the image of an infinite rich word $\uu \in \Sigma^\omega_d$ under the morphism~$\eta  : i\mapsto r_i$. This morphism belongs to  the class $P_{\mathrm{ret}}$ and $p$ is its marker. 
\end{itemize}

\section{The Compatibility Lemma and Compatible Sets}\label{Section:UsefulLemma}

Recall that by \autoref{Theorem:CompleteReturnsToPalindromes}, every complete return word to a palindrome in an infinite rich word is itself a palindrome.  The following lemma says even more about the structure of return words to palindromes in recurrent rich words, and will be very useful for establishing both the ``outer layer'' and the ``inner layers'' of the structure theorems that we are concerned with.

\begin{lemma}[The Compatibility Lemma]
\makeatletter\def\@currentlabelname{Compatibility Lemma}\makeatother
\label{Lemma:Compatibility}  
Let $\mathbf{w} \in \Sigma_d^\omega$ be a recurrent rich word and   $q \in \Sigma_d^
*$ be a palindrome.  If  $rq$ and $r'q$  are two complete return words to $q$ in $\mathbf{w}$ containing   a palindrome   $s\in \Sigma_d^*$, then the longest common prefix of $rq$  and $r'q$ contains $s$. 
\end{lemma}

\begin{proof} Assume that $s$ does not occur in $q$, otherwise the statement is trivial.

Let $r_0r_1\cdots r_{n}$ be a factor of $\mathbf{w}$  composed of return words $r_0, r_1, \ldots, r_{n}$  to $q$ in $\mathbf{w}$, where $n \in \N, n>0$.   Assume that $s$ occurs only in the complete return words $r_0q$  and  $r_{n}q$ but not in $r_iq$ for $i=1,2, \ldots, n-1$.  

Let $P$ be the shortest suffix of $r_0q$ containing  $s$
and let $P'$ be the shortest prefix of $r_nq$  containing  $s$.  
As $s$ does not occur in $q$, we deduce that $q$ is a proper suffix of $P$ and a proper prefix of $P'$. 

Let $v$ be a  factor of $r_0r_1 \cdots r_nq$ such that $P$ is its prefix and $P'$ is its suffix. In $v$,  there are exactly two occurrences of the palindrome  $s$,   one at the beginning of the word $v$ and one at its  end. In other words,  $v$ is a complete return word to $s$ in $\ww$.  By \autoref{Theorem:CompleteReturnsToPalindromes}, it follows that $v$ is a palindrome, and thus $P' = R(P)$. 
  
Recall that the complete return word $r_0q$ is a palindrome  with a suffix $P$.  Hence,  $R(P)$ is a prefix of $r_0q$.   So $r_0q$ and $r_nq$ share a common prefix $P'$.  

To sum it up, we have shown that the two nearest return words containing the factor $s$ share a common prefix $P'$ which contains $s$. Since this holds true for every pair of consecutive complete return  words containing $s$, the lemma is proven.
\end{proof}

The \nameref{Lemma:Compatibility} motivates the following terminology.

\begin{definition}\label{Definition:Compatible}
    For a palindrome $q$ and two return words $r$ and $r'$ to $q$, we say that $r$ and $r'$ are \emph{compatible} if for every palindrome $s$ such that $s$ occurs in both $rq$ and $r'q$, the longest common prefix of $rq$ and $r'q$ contains $s$.

    We say that a set of return words to a palindrome $q$ is \emph{compatible} if its elements are pairwise compatible.
\end{definition}

With this new terminology, we can restate the \nameref{Lemma:Compatibility} as follows.

\begin{observation}\label{Observation:CompatibleSets}
    If $q$ is a palindromic factor of a recurrent rich word $\zz$, then $\Ret_\zz(q)$ is compatible.
\end{observation}

The next three examples demonstrate techniques that will be vital to our proofs of the structural parts of \autoref{Theorem:TernaryStructure} and \autoref{Theorem:QuaternaryStructure}.  For the benefit of the reader, the first two examples include illustrations and can be verified by hand.  The third example illustrates the computer-assisted approach that we will take to establish similar statements in the remainder of the paper.

\begin{example}
    Let $\zz\in\Sigma_3^\omega$ be a recurrent rich word containing $\tt{0}$, and consider return words to $\tt{0}$ in $\zz$.  The words $\tt{01}$ and $\tt{0212}$ are not compatible, as the palindrome $\tt{1}$ occurs in both $\tt{010}$ and $\tt{02120}$, but the longest common prefix of $\tt{010}$ and $\tt{02120}$ is $\tt{0}$, which does not contain $\tt{1}$.  On the other hand, it is easy to check that $\tt{022}$ and $\tt{0212}$ are compatible.

    \begin{figure}
        \centering
        \begin{tikzpicture}[scale=0.8]
            \node (empty) at (0,0) {\scriptsize $\emptyset$};
            \node (01) at (-1,2) {\scriptsize $\{\tt{01}\}$};
            \node (0121) at (-3,2) {\scriptsize $\{\tt{0121}\}$};
            \node (012121) at (-5,2) {\scriptsize $\{\tt{012121}\}$};
            \node (02) at (1,2) {\scriptsize $\{\tt{02}\}$};
            \node (022) at (3,2) {\scriptsize $\{\tt{022}\}$};
            \node (0212) at (5,2) {\scriptsize $\{\tt{0212}\}$};
            \node (01+0121) at (-2,4) {\scriptsize $\{\tt{01,0121}\}$};
            \node (01+012121) at (-4.5,4) {\scriptsize $\{\tt{01,012121}\}$};
            \node (0121+012121) at (-7.5,4) {\scriptsize $\{\tt{0121,012121}\}$};
            \node (01+02) at (0,4) {\scriptsize $\{\tt{01,02}\}$};
            \node (01+022) at (2,4) {\scriptsize $\{\tt{01,022}\}$};
            \node (02+022) at (4,4) {\scriptsize $\{\tt{02,022}\}$};
            \node (02+0212) at (6,4) {\scriptsize $\{\tt{02,0212}\}$};
            \node (022+0212) at (8,4) {\scriptsize $\{\tt{022,0212}\}$};
            \node (01+0121+012121) at (-4.5,6) {\scriptsize $\{\tt{01,0121,012121}\}$};
            \node (01+02+022) at (2,6) {\scriptsize $\{\tt{01,02,022}\}$};
            \node (02+022+0212) at (6,6) {\scriptsize $\{\tt{02,022,0212}\}$};
            \draw (empty) -- (01);
            \draw (empty) -- (0121);
            \draw (empty) -- (012121);
            \draw (empty) -- (02);
            \draw (empty) -- (022);
            \draw (empty) -- (0212);
            \draw (012121) -- (0121+012121);
            \draw (012121) -- (01+012121);
            \draw (0121) -- (0121+012121);
            \draw (0121) -- (01+0121);
            \draw (01) -- (01+012121);
            \draw (01) -- (01+0121);
            \draw (01) -- (01+02);
            \draw (01) -- (01+022);
            \draw (02) -- (01+02);
            \draw (02) -- (02+022);
            \draw (02) -- (02+0212);
            \draw (022) -- (01+022);
            \draw (022) -- (02+022);
            \draw (022) -- (022+0212);
            \draw (0212) -- (02+0212);
            \draw (0212) -- (022+0212);
            \draw (01+0121+012121) -- (01+0121);
            \draw (01+0121+012121) -- (01+012121);
            \draw (01+0121+012121) -- (0121+012121);
            \draw (01+02+022) -- (01+02);
            \draw (01+02+022) -- (01+022);
            \draw (01+02+022) -- (02+022);
            \draw (02+022+0212) -- (02+022);
            \draw (02+022+0212) -- (02+0212);
            \draw (02+022+0212) -- (022+0212);
        \end{tikzpicture}
        \caption{The compatible subsets of $S=\{\tt{01},\tt{0121},\tt{012121},\tt{02},\tt{0212},\tt{022}\}$, ordered by inclusion.}
        \label{Figure:MaximalCompatibleSets}
    \end{figure}
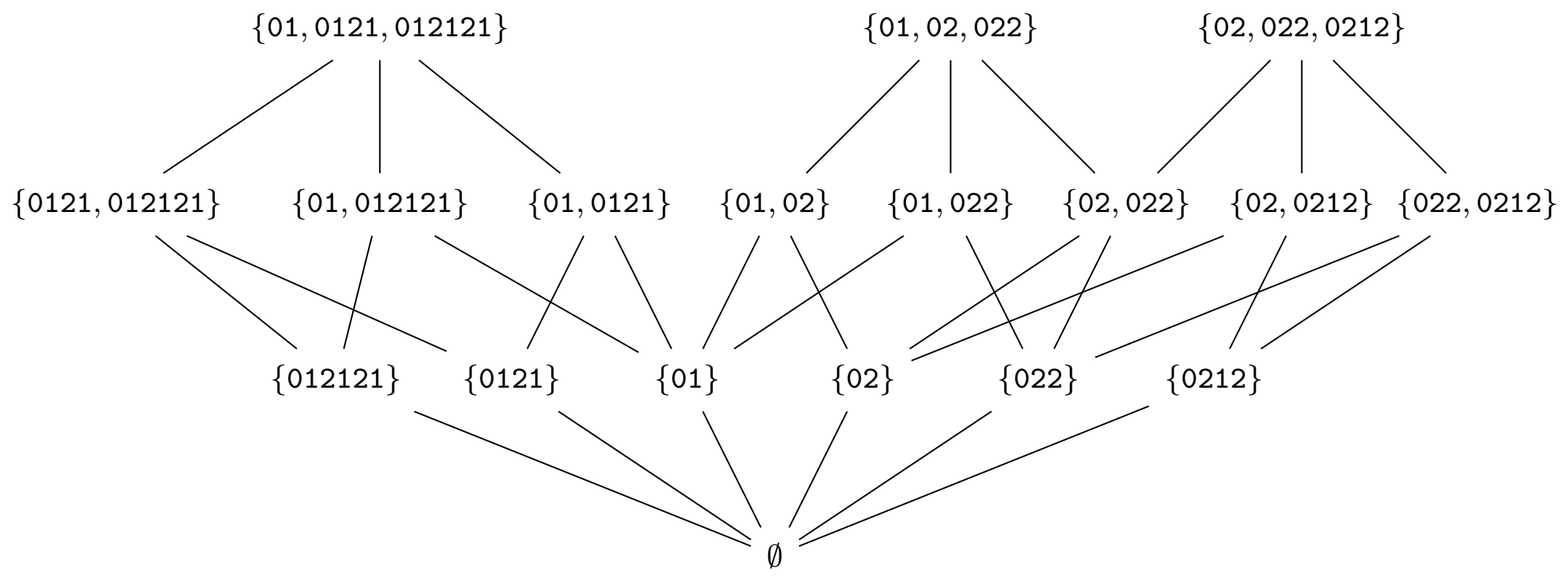

    Suppose now that we are told that $\mathcal{R}_\mathbf{z}(\tt{0})$ is a subset of
    \[
        S=\{\tt{01},\tt{0121},\tt{012121},\tt{02},\tt{0212},\tt{022}\}.
    \]
    By \autoref{Observation:CompatibleSets}, we know that $\mathcal{R}_\zz(\tt{0})$ must be a compatible subset of $S$.  The compatible subsets of $S$ form a poset with respect to inclusion, whose Hasse diagram is illustrated in \autoref{Figure:MaximalCompatibleSets}.  We may conclude that $\mathcal{R}_\zz(\tt{0})$ is contained in one of the maximal compatible subsets of $S$, which are $\{\tt{01},\tt{0121},\tt{012121}\}$, $\{\tt{01},\tt{02},\tt{022}\}$, and $\{\tt{02},\tt{022},\tt{0212}\}$. 
\end{example}

\begin{example}
    Let $\zz$ be a recurrent $\tfrac{14}{5}$-power-free binary rich word, and suppose that $\zz$ contains the factor $w=\tt{0011011}$.  We will use the \nameref{Lemma:Compatibility} to find a set that contains all complete return words to $\tt{00}$ in $\zz$.  First, we find the set $S$ of all palindromic factors of $w$.  We have
    \[
        S=\{\tt{0},\tt{00},\tt{1},\tt{11},\tt{0110},\tt{101},\tt{11011}\}.
    \]
    For each $s\in S$, let $p_s$ denote the shortest prefix of $w$ containing $s$, e.g., we have $p_\tt{1}=\tt{001}$ and $p_\tt{101}=\tt{001101}$.  By the \nameref{Lemma:Compatibility}, we know that every complete return word to $\tt{00}$ in $\zz$ containing $s\in S$ must have prefix $p_s$.

    To find all possible complete return words to $\tt{00}$ in $\zz$, we perform a depth-first search through the tree of words over $\Sigma_2$ with prefix $\tt{00}$.  Each time that we add a letter, we check that the resulting word $u$ is rich and $\tfrac{14}{5}$-power-free, and that if $u$ has suffix $s\in S$, then $u$ has prefix $p_s$.  If at least one of these conditions is not satisfied, then we backtrack.  If all of these conditions are satisfied and $u$ has suffix $\tt{00}$, then we have a possible complete return word to $\tt{00}$ in $\zz$, so we add $u$ to a set $X$ (which was initially empty), and then backtrack.  

        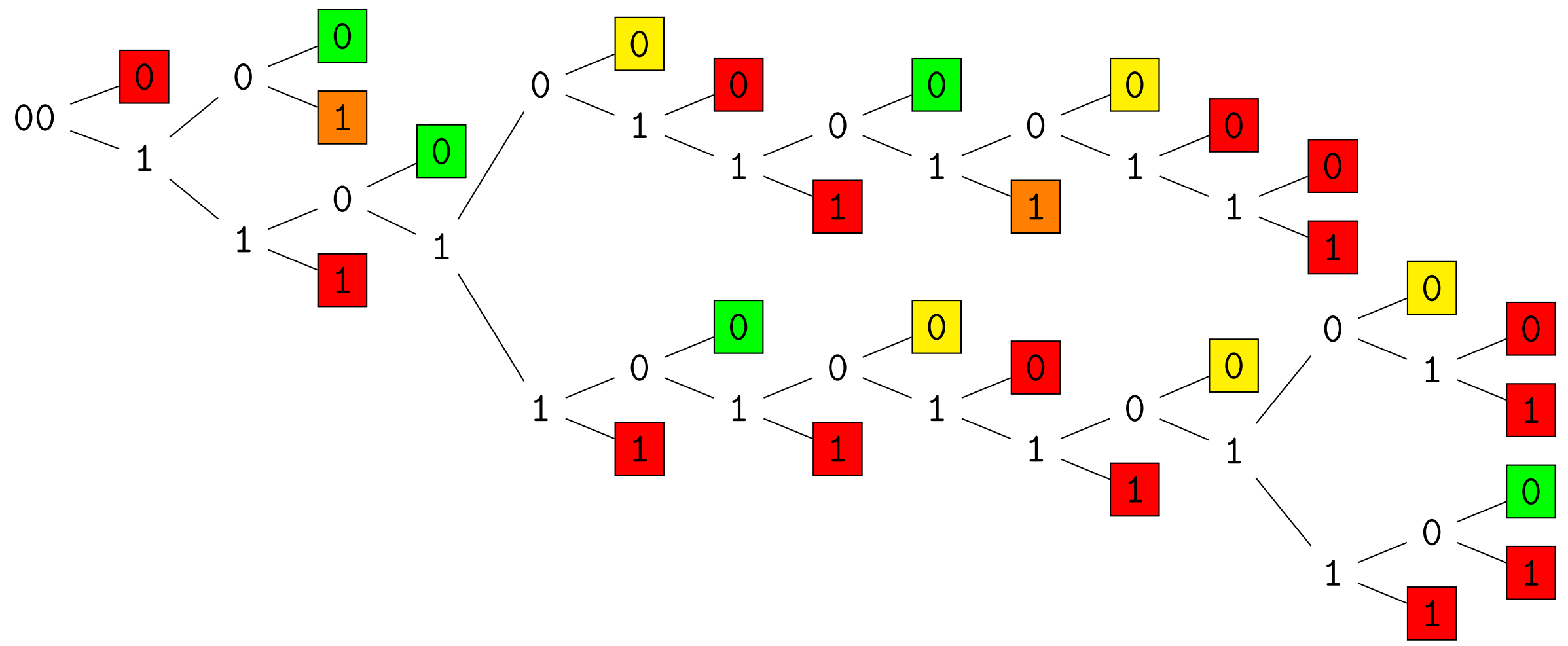
\begin{figure}
        \centering
        \begin{forest}
            [
                $\tt{00}$,for tree={grow=0,l=1cm}
                [
                    $\tt{1}$
                    [
                        $\tt{1}$
                        [
                            $\tt{1}$,draw,fill=red
                        ]
                        [
                            $\tt{0}$
                            [
                                $\tt{1}$
                                [
                                    $\tt{1}$
                                    [
                                        $\tt{1}$,draw,fill=red
                                    ]
                                    [
                                        $\tt{0}$
                                        [
                                            $\tt{1}$
                                            [
                                                $\tt{1}$,draw,fill=red
                                            ]
                                            [
                                                $\tt{0}$
                                                [
                                                    $\tt{1}$
                                                    [
                                                        $\tt{1}$
                                                        [
                                                            $\tt{1}$,draw,fill=red
                                                        ]
                                                        [
                                                            $\tt{0}$
                                                            [
                                                                $\tt{1}$
                                                                [
                                                                    $\tt{1}$
                                                                    [
                                                                        $\tt{1}$,draw,fill=red
                                                                    ]
                                                                    [
                                                                        $\tt{0}$
                                                                        [   
                                                                            $\tt{1}$,draw,fill=red
                                                                        ]
                                                                        [
                                                                            $\tt{0}$,draw,fill=green
                                                                        ]
                                                                    ]
                                                                ]
                                                                [
                                                                    $\tt{0}$
                                                                    [
                                                                        $\tt{1}$
                                                                        [
                                                                            $\tt{1}$,draw,fill=red
                                                                        ]
                                                                        [
                                                                            $\tt{0}$,draw,fill=red
                                                                        ]
                                                                    ]
                                                                    [
                                                                        $\tt{0}$,draw,fill=yellow
                                                                    ]
                                                                ]
                                                            ]
                                                            [
                                                                $\tt{0}$,draw,fill=yellow
                                                            ]
                                                        ]
                                                    ]
                                                    [
                                                        $\tt{0}$,draw,fill=red
                                                    ]
                                                ]
                                                [
                                                    $\tt{0}$,draw,fill=yellow
                                                ]
                                            ]
                                        ]
                                        [
                                            $\tt{0}$,draw,fill=green
                                        ]
                                    ]
                                ]
                                [
                                    $\tt{0}$
                                    [
                                        $\tt{1}$
                                        [
                                            $\tt{1}$
                                            [
                                                $\tt{1}$,draw,fill=red
                                            ]
                                            [
                                                $\tt{0}$
                                                [
                                                    $\tt{1}$
                                                    [
                                                        $\tt{1}$,draw,fill=orange
                                                    ]
                                                    [
                                                        $\tt{0}$
                                                        [
                                                            $\tt{1}$
                                                            [
                                                                $\tt{1}$
                                                                [
                                                                    $\tt{1}$,draw,fill=red
                                                                ]
                                                                [
                                                                    $\tt{0}$,draw,fill=red
                                                                ]
                                                            ]
                                                            [
                                                                $\tt{0}$,draw,fill=red
                                                            ]
                                                        ]
                                                        [
                                                            $\tt{0}$,draw,fill=yellow
                                                        ]
                                                    ]
                                                ]
                                                [
                                                    $\tt{0}$,draw,fill=green
                                                ]
                                            ]
                                        ]
                                        [
                                            $\tt{0}$,draw,fill=red
                                        ]
                                    ]
                                    [
                                        $\tt{0}$,draw,fill=yellow
                                    ]
                                ]
                            ]
                            [
                                $\tt{0}$,draw,fill=green
                            ]
                        ]
                    ]
                    [
                        $\tt{0}$
                        [
                            $\tt{1}$,draw,fill=orange
                        ]
                        [
                            $\tt{0}$,draw,fill=green
                        ]
                    ]
                ]
                [
                    $\tt{0}$,draw,fill=red
                ]
            ]
        \end{forest}
        \caption{A backtracking search for complete return words to $\tt{00}$ compatible with $\tt{0011011}$ in a recurrent $\tfrac{14}{5}$-power-free rich word $\zz$.}
        \label{Figure:BinaryReturns}
    \end{figure}%
    
    The search is illustrated with a tree in \autoref{Figure:BinaryReturns}.  The colour of each leaf indicates the reason for backtracking, as described below.
    \begin{itemize}
        \item Each green leaf corresponds to a possible complete return word to $\tt{00}$ in $\zz$.
        \item Each red leaf corresponds to a word with a suffix of exponent at least $\tfrac{14}{5}$.
        \item Each yellow leaf corresponds to a word that is not rich.
        \item Each orange leaf corresponds to a word with suffix $s$ for some $s\in S$ that does not have prefix $p_s$.  In particular, the word $\tt{00101}$ contains $\tt{101}\in S$ but does not have prefix $p_\tt{101}=\tt{001101}$, while the word $\tt{001101011011}$ contains $\tt{11011}\in S$ but does not have prefix $p_{\tt{11011}}=\tt{0011011}$.
    \end{itemize}
    Thus, the search terminates with
    \[
        X=\{\tt{00100},\tt{001100},\tt{00110101100},\tt{001101100},\tt{00110110101101100}\},
    \]
    and we conclude that 
    \[
        \mathcal{R}_\zz(\tt{00})\subseteq \{\tt{001},\tt{0011},\tt{001101011},\tt{0011011},\tt{001101101011011}\}. \qedhere
    \]
\end{example} 

We call the depth-first search for possible complete return words to $\tt{00}$ in $\zz$ used in \autoref{Example:OuterLayerForTernaryRichWords} a ``backtracking search for complete return words to $\tt{00}$ in $\zz$ compatible with the factor $w$''. 
We can perform a similar search for possible complete return words to any palindrome $p$ in any recurrent rich word $\xx$ if we know a prefix of a complete return word to $p$ in $\xx$.
Note that such searches are not necessarily guaranteed to terminate in general; a restrictive condition on the critical exponent of $\xx$ is helpful in this regard.  
Importantly, it is straightforward to automate these searches and have them performed entirely by computer, as in the next example, where we establish a (weak) ``outer layer'' for \autoref{Theorem:TernaryStructure}.

\begin{example}\label{Example:OuterLayerForTernaryRichWords}
    Let $\zz$ be a recurrent $\tfrac{16}{7}$-power-free rich word over $\Sigma_3$.  A backtracking search shows that a longest $\tfrac{16}{7}$-power-free rich word over $\Sigma_3$ containing no factor obtained by permuting the letters of 
    \[
        w=\tt{001011010020220200202202001002}
    \]
    has length $681$.  So we may assume, without loss of generality, that $\zz$ contains the word~$w$. 
    
    By a backtracking search for complete return words to $p=\tt{0010110100}$ in $\zz$ compatible with $w$, we find that 
    \[
        \Ret_\mathbf{z}(q)\subseteq \{A,AD,AC,ACC,ACCBCC\}.
    \]
    where $A=\tt{00101101}$, $B=\tt{001}$, $C=\tt{00202202}$, and $D=\tt{002}$, as defined in \autoref{Theorem:TernaryStructure}.  Thus, up to permutation of the letters, some suffix of $\zz$ must have the form $\pi(\mathbf{x})$ for some infinite word $\mathbf{x}\in \Sigma_5^\omega$, where $\pi$ is the morphism defined in \autoref{Theorem:TernaryStructure}.  
\end{example}

\section{Persistent Repetitions}\label{Section:PersistentRepetitions}

For $t\geq 1$ and $x\in\Sigma_d^*$ with $x\neq \varepsilon$, we say that a word $y\in\Sigma_d^*$ is a \emph{$t$-power of $x$} if $y$ is a prefix of $x^\omega$ and $|y|\geq t |x|$.  We say that $y$ is a \emph{$t$-power} if it is a $t$-power of some word $x$. Note that this definition of $t$-power is non-standard. Under our definition, the word $y$ is a $t$-power if and only if it has exponent \emph{at least} $t$.  Some authors say that $y$ is a $t$-power if and only if it has exponent \emph{exactly} $t$.

Throughout this section, we let  
$\varphi=\varphi_d$ for some fixed $d\geq 2$, and we take $\rho: \Sigma_d^* \mapsto \Sigma_{d'}^*$ to be a morphism satisfying  
 \begin{equation}\label{eq:increase}|\rho(\tt{0})|\leq |\rho(\tt{1})|\leq \cdots\leq |\rho(d-1)|.\end{equation}
Our goal in this section is to identify factors whose occurrence in an infinite word  $\bf {x}$ causes the occurrence of a $t$-power in the infinite word $\rho(\varphi^n(\bf {x}$)) for every $n \in \mathbb{N}$ and  suitably chosen $t> E^*(\uu_d)$.  Eliminating such factors will facilitate our proof of the ``inner layers'' of our structure theorems.

\subsection{Bookended repetitions}

First, we study how powers occurring in the language of an infinite word behave after applying a morphism. 
From the context in which the power  is situated, we can sometimes also deduce the context of its image. By the context of a factor we mean the letter preceding the factor and the letter following it. However, the context is not always known.

The following observation is inspired by the technique of Klouda \cite{K2012} for finding all bispecial factors in a fixed point of a morphism.  
For a set of words $X$, we let $\lcp(X)$ and $\lcs(X)$ denote the longest common prefix and longest common suffix, respectively, of all words in $X$.

\begin{observation}\label{lem:new} Let $\eta: \Sigma^*_d \mapsto \Sigma_{d'}^*$ be a morphism,  let $g,h \in \Sigma_d \cup\{\varepsilon\}$, and let $z$ be a finite prefix of $x^\omega$ of length $\geq|x|$ for some nonempty $x \in \Sigma_d^*$, i.e.,  
$z=xy=yx'$ for some words $x',y\in\Sigma_d^*$ with $|x|=|x'|$. 
Let $g'$ denote the last letter of $x$ and $h'$ denote the first letter of $x'$, and let $s =\lcs\{\eta(c): c \in \Sigma_d\}$ and $p =\lcp\{\eta(c): c \in \Sigma_d\}$.  

\medskip

\noindent  Suppose that $gzh$ occurs in  a recurrent word ${\bf x} \in \Sigma_d^\omega$.  Then 
\begin{enumerate}  \item $S\eta(z)P$  occurs in $\eta({\bf x})$, where  $ S=   {\lcs}\{s\eta(g), s\eta(g')\}$ and $   P =  {\lcp}\{\eta(h)p, \eta(h')p\}$;

\item  $S\eta(z)P$  is a prefix of  $\bigl(S\eta(x)S^{-1}\bigr)^\omega$.
\end{enumerate}
Moreover, if $|s\eta(g)| > |S|$,  let $G$ denote the letter such that  $GS$ is a suffix of $s\eta(g)$, and let $G = \varepsilon$ otherwise. Analogously,  if  $|\eta(h)p| > |P|$, let $H$ denote the letter such that $PH$ is a~prefix of $\eta(h)p$, and let $H = \varepsilon$ otherwise.  Then $GS\eta(z)PH$ occurs in $\eta({\bf x})$. 

\end{observation}

\begin{example}\label{Example:BookendedRepetitions}
    Let $\varphi=\varphi_3$, and let $z=\tt{020}$, $x=\tt{02}$, and $g=h=\tt{1}$.  Using the notation of \autoref{lem:new}, we have $s=\varepsilon$, $p=\tt{0}$, $g'=\tt{2}$, and $h'=\tt{2}$.  Thus we have $S=\lcs\{\tt{02},\tt{022}\}=\tt{2}$ and $P=\lcp\{\tt{020},\tt{0220}\}=\tt{02}$.  Thus, if a recurrent word $\xx$ contains $gzh=\tt{10201}$, then $\varphi(\xx)$ contains $S\varphi(z)P=\tt{2010220102}$.  We deduce in particular that if $\varphi(\xx)$ is square-free, then $\xx$ does not contain the factor $\tt{10201}$.  Finally, we have $G=\tt{0}$ and $H=\tt{0}$, so in fact, if $\xx$ contains $gzh=\tt{10201}$, then $\varphi(\xx)$ contains $GS\pi(z)PH=\tt{020102201020}$, i.e., we know the context in which the repetition $S\pi(z)P=(\tt{20102})^2$ occurs in $\varphi(\xx)$.

    As another example (still with $\varphi=\varphi_3$), let $z=\tt{00}$, $x=\tt{0}$, and $g=h=\varepsilon$.  Then we have $S=\lcs\{\varepsilon,\tt{01}\}=\varepsilon$ and $P=\lcp\{\tt{0},\tt{010}\}=\tt{0}$.  Thus, if $\xx$ contains the factor $gzh=\tt{00}$, then $\varphi(\xx)$ contains the factor $S\varphi(z)P=\tt{01010}$.  In this case, however, we do not know the context in which the repetition $\tt{01010}$ occurs in $\varphi(\xx)$, because we have $G=H=\varepsilon$.
\end{example}


Let us introduce some terminology and notation.  The quadruple $R = (g,z,h,x)$  of \autoref{lem:new} is called a \emph{bookended repetition} over $\Sigma_d$.  The quadruple $(G,S\eta(z)P,H, S\eta(x)S^{-1})$ is also a~bookended repetition over $\Sigma_{d'}$, and we will abuse notation slightly and denote it simply by $\eta(R)$.  With this new notation, in \autoref{Example:BookendedRepetitions}, we found that
\[
    \varphi(\tt{1},\tt{020},\tt{1},\tt{02})=(\tt{0},\tt{2010220102},\tt{0},\tt{20102})
\]
and
\[
    \varphi(\varepsilon,\tt{00},\varepsilon,\tt{0})=(\varepsilon,\tt{01010},\varepsilon,\tt{01}).
\]
We say that a bookended repetition $R=(g,z,h,x)$ is a~\emph{$t$-power} if $z$ is a $t$-power, and we say that a word $\xx$ \emph{contains} $R$ if the word $gzh$ is a factor of $\xx$.  So from \autoref{Example:BookendedRepetitions}, we see that if $\xx$ contains the bookended repetition $R=(\tt{1},\tt{020},\tt{1},\tt{02})$, then $\varphi(\xx)$ contains the $2$-power $\varphi(R)=(\tt{0},\tt{2010220102},\tt{0},\tt{20102})$.

\begin{definition} Let  $t>1$ and $\rho: \Sigma_d\mapsto \Sigma_{d'}$ satisfy \eqref{eq:increase}. 
    We say that a bookended repetition $R$ over $\Sigma_d$ is \emph{$t$-persistent (with respect to $\varphi$ and $\rho$)} if $\rho(\varphi^n(R))$ is a $t$-power for all $n\geq 0$.
\end{definition}

\begin{observation}\label{Observation:ProtoPowers}
    Let $\zz$ be a recurrent word of the form $\zz=\rho(\varphi^n(\yy))$ for some $n\geq 0$ and $\yy\in\Sigma_d^*$.  If $\zz$ is $t$-power-free, then $\yy$ contains no $t$-persistent bookended repetition.
\end{observation}

While it is straightforward to determine whether or not a single given bookended repetition is a~$t$-power, deciding whether or not all bookended repetitions in the infinite sequence $\rho(R), \rho(\varphi(R)), \rho(\varphi^2(R)), \ldots$ are $t$-powers (i.e., deciding whether or not $R$ is $t$-persistent) may not be easy, or even possible.  Verifying that $R$ is $t$-persistent will be facilitated if a~special type of bookended repetition that we introduce below appears in the sequence  $R$, $\varphi(R)$, $\varphi^2(R), \ldots$.

\subsection{Top-heavy \texorpdfstring{$t$}{t}-powers}

\begin{definition}
Let $t>1$ and $y$ be a finite prefix of $x^\omega$ for some nonempty $x \in \Sigma_d^*$. We say that $y$ is a~\emph{top-heavy $t$-power} of $x$  if all of the following inequalities concerning the number of occurrences of letters in $x$ and $y$ are satisfied.   
    \begin{align*}
        |y|_{d-1}&\geq t|x|_{d-1}\\
        |y|_{d-2}+|y|_{d-1}&\geq t\left(|x|_{d-2}+|x|_{d-1}\right)\\
         & \ \, \vdots\\
        |y|_\tt{0}+|y|_\tt{1}+\cdots+|y|_{d-1}&\geq t\left(|x|_\tt{0}+|x|_\tt{1}+\cdots +|x|_{d-1}\right)
    \end{align*}
    In particular, the last inequality says that $|y|\geq t|x|$, from which it follows that if $y$ is a~top-heavy $t$-power of $x$, then $y$ is an (ordinary) $t$-power of $x$.
\end{definition}

For example, consider the word $w=\tt{010201}=(\tt{0102})^{3/2}=x^{3/2}\in\Sigma_3^*$.  Note that $w$ is a $\tfrac{3}{2}$-power of $x=\tt{0102}$, but not a top-heavy $\tfrac{3}{2}$-power of $x$, because $|w|_\tt{2}\not\geq \tfrac{3}{2}|x|_\tt{2}$.  
Now consider the word $W=\varphi_3(w)$.  We have
\[
W=\tt{0102010220102}=(\tt{010201022})^{13/9}=X^{13/9},
\]
where $X=\varphi(x)=\tt{010201022}$.  Again, we see that $W$ is a $\tfrac{13}{9}$-power of $X$, but not a~top-heavy $\tfrac{13}{9}$-power of $X$.  Note, however, that $W$ is a top-heavy $\tfrac{4}{3}$-power of $X$, since 
\begin{align*}
    |W|_\tt{2}&=4=\tfrac{4}{3}|X|_2\,,\\
    |W|_\tt{1}+|W|_\tt{2}&=3+4=7>\tfrac{20}{3}=\tfrac{4}{3}(|X|_\tt{1}+|X|_{\tt{2}})\,, \text{ and }\\
    |W|_\tt{0}+|W|_\tt{1}+|W|_\tt{2}&=|W|=13> 12=\tfrac{4}{3}|X|=\tfrac{4}{3}(|X|_\tt{0}+|X|_\tt{1}+|X|_\tt{2})\,.
\end{align*}
The terminology ``top-heavy'' suggests the essential intuition here, which is that the ``importance/weight'' associated with each letter increases with the size of the letter.  

The next lemma says that the image of any top-heavy $t$-power under a morphism satisfying (\ref{eq:increase}) must be a $t$-power.

\begin{lemma}\label{Lemma:TopHeavyImages}
    Let $\eta:\Sigma_d^*\rightarrow \Sigma_{d'}^*$ be a morphism satisfying $|\eta(\tt{0})|\leq |\eta(\tt{1})|\leq \cdots\leq |\eta(d-1)|$.
    If $y\in \Sigma_d^*$ is a top-heavy $t$-power of some nonempty word $x \in \Sigma_d^*$, then $\eta(y)$ is a $t$-power of $\eta(x)$.
\end{lemma}

In order to prove \autoref{Lemma:TopHeavyImages}, we make use of a simple claim. It is based on a summation method due to Abel.

\begin{claim}[Abel's summation]\label{claim:Abel} 
Let $0\leq\alpha_0\leq \alpha_1 \leq \cdots \leq \alpha_{d-1}$ and $\beta_0, \beta_1, \ldots, \beta_{d-1} $ be real numbers such that $B_i := \displaystyle\sum_{k=i}^{d-1}\beta_k\geq 0$ for every $i =0,1,\ldots, d$. Then $\displaystyle\sum_{k=0}^{d-1} \alpha_k\beta_k \geq 0$.    
\end{claim}
\begin{proof} 
We have
\[
\sum_{k=0}^{d-1} \alpha_k\beta_k =  \sum_{k=0}^{d-1} \alpha_k(B_k-B_{k+1})= \sum_{k=0}^{d-1} \alpha_kB_k - \sum_{k=1}^{d-1} \alpha_{k-1}B_k\,,
\]
as $B_d$ equals $0$ (the empty sum). Hence 
\[
\sum_{k=0}^{d-1} \alpha_k\beta_k = \alpha_0 B_0 + \sum_{k=1}^{d-1} \underbrace{(\alpha_k - \alpha_{k-1})}_{\geq 0}B_k \geq 0\,.   \qedhere
\]
\end{proof}

\begin{proof}[Proof of \autoref{Lemma:TopHeavyImages}]
As $y$ is a prefix of $x^\omega$, obviously, $\eta(y)$ is a prefix  of $\bigl(\eta(x)\bigr)^\omega$. 
For each $k\in\Sigma_d$,  we let $\alpha_k =|\eta(k)|$. As $|\eta(y)|= \displaystyle\sum_{k=0}^{d-1}\alpha_k|y|_k$ and, analogously, $|\eta(x)|= \displaystyle\sum_{k=0}^{d-1}\alpha_k|x|_k$, we need to show that 
\begin{align}\label{AbelSumInequality}
 \sum_{k=0}^{d-1}\alpha_k\bigl(|y|_k-t|x|_k\bigr)\geq 0
\,.
\end{align}
Let $\beta_k = |y|_k-t|x|_k$, and observe that for each $i\in\{0,1,\ldots,d-1\}$, we have
\[
\sum_{k=i}^{d-1}\beta_i=\sum_{k=i}^{d-1}|y|_k-t\sum_{k=i}^{d-1}|x|_k\geq 0\,,
\]
since $y$ is a top heavy $t$-power of $x$.  Therefore, the $\alpha_k$ and $\beta_k$ satisfy the assumptions of  \autoref{claim:Abel}, from which inequality (\ref{AbelSumInequality}) follows immediately.  
\end{proof}

\begin{proposition}[The Persistence Criterion]
\makeatletter\def\@currentlabelname{Persistence Criterion}\makeatother
\label{Proposition:Persistence}  Let  $t> 1$ and $\rho: \Sigma_d^* \mapsto \Sigma_{d'}^*$  be a morphism satisfying  \eqref{eq:increase}.  Assume that $R$ is a bookended repetition over $\Sigma_d$.   If there exists $\ell \in \mathbb{N}$ such that $\rho(\varphi^i(R))$ is a~$t$-power for every $i=0,1,\ldots, \ell-1$  and $\varphi^\ell(R)$ is a~top-heavy $t$-power, then  $R$
 is $t$-persistent. 
\end{proposition}

\begin{proof}
    We must show that for all $n\geq 0$, the bookended repetition 
    \[
        \rho(\varphi^{\ell+n}(R))=\rho(\varphi^n(\varphi^\ell(R)))
    \]
    is a $t$-power.  By \autoref{Lemma:TopHeavyImages}, since $\varphi^\ell(R)$ is a top-heavy $t$-power,
    it suffices to show that for all $n\geq 0$, we have 
   \begin{equation}\label{eq:composition}
    |\rho(\varphi^n(\tt{0}))|\leq |\rho(\varphi^n(\tt{1}))|\leq \cdots\leq |\rho(\varphi^n(d-2))|\leq |\rho(\varphi^n(d-1))|.
\end{equation}

    We prove this statement by induction  on $n$ for every $n \geq 0$.  The base case $n=0$ is immediate.
    Now suppose for some integer $n\geq 0$ that the inequalities \eqref{eq:composition} take place. 
    For all $i\in\Sigma_d$ with $i< d-2$, we have 
        \begin{align*}
            |\rho(\varphi^{n+1}(i))|
            &=|\rho(\varphi^n(\tt{0}))|+|\rho(\varphi^n(i+1))|\\
            &\leq |\rho(\varphi^n(\tt{0}))|+|\rho(\varphi^n(i+2))|\\
            &=|\rho(\varphi^{n+1}(i+1))|.
        \end{align*}
As $\varphi(d-2)$ is a proper prefix of $\varphi(d-1)$, the inequality $|\rho(\varphi^{n+1}(d-2))| \leq |\rho(\varphi^{n+1}(d-1))|$ is obvious. 
\end{proof}

The \nameref{Proposition:Persistence} essentially tells us that it is easy to decide whether or not a bookended repetition $R$ is $t$-persistent if we can find a~top heavy $t$-power in the sequence $R, \varphi(R), \varphi^2(R), \ldots$.   The minimum $\ell\in \mathbb{N}$ such that $\varphi^\ell(R)$  is a top heavy $t$-power (if such an $\ell$ exists) is called the \emph{tipping point} of $R$.  For a fixed integer $k$, a procedure for deciding whether or not a given bookended repetition $R$ is $t$-persistent with tipping point at most $k$ is given explicitly in \autoref{Algorithm:EventuallyTopHeavy}.  Restricting the tipping point to be at most $k$ simply ensures that we do not get stuck in an infinite loop.  This procedure will be used extensively in the proof of the structural parts of \autoref{Theorem:TernaryStructure} and \autoref{Theorem:QuaternaryStructure}.  In our implementation of this procedure, we work exclusively with the Parikh vectors of the words involved, which simplifies calculations greatly.

\begin{algorithm} 
\caption{Deciding whether or not a bookended repetition $R$ is $t$-persistent (with respect to $\varphi$ and $\rho$) with tipping point at most $k$}\label{Algorithm:EventuallyTopHeavy} 
\begin{algorithmic}[1] 
\Procedure{IsPersistent}{$R,t,\varphi,\rho,k$} 
\State $\ell=0$ 
\While{$\ell\leq k$}
\If{$R$ is a top-heavy $t$-power}
    \State \textbf{return} \texttt{True}
\EndIf
\If{$\rho(R)$ is not a $t$-power}
    \State \textbf{return} \texttt{False}
\EndIf
\State $R\gets \varphi(R)$ 
\State $\ell\gets \ell+1$
\EndWhile
\State \textbf{return} \texttt{False} 
\EndProcedure 
\end{algorithmic} 
\end{algorithm}

We conclude this section by demonstrating the procedure of \autoref{Algorithm:EventuallyTopHeavy} on a concrete example which is small enough to complete by hand.

\begin{example}
    Let $\varphi=\varphi_3$ and $\rho=\pi_1$, where $\pi_1$ is the morphism defined in \autoref{Theorem:BinaryStructure}, i.e., we have $\rho(\tt{0})=\tt{0}$, $\rho(\tt{1})=\tt{01}$, and $\rho(\tt{2})=\tt{011}$.  Consider the bookended repetition $R=(\tt{2},\tt{010},\tt{2},\tt{01})$.  We will show that $R$ is $\tfrac{7}{3}$-persistent with respect to $\varphi$ and $\rho$ with tipping point $2$.

    First note that while $R$ itself is not a $\tfrac{7}{3}$-power, its image $\rho(R)=(\tt{1},\tt{1001001},\tt{1},\tt{100})$ is a $\tfrac{7}{3}$-power.  Similarly, while $\varphi(R)=(\tt{2},\tt{201020102},\tt{2},\tt{2010})$
    is not a $\tfrac{7}{3}$-power, its image $\rho(\varphi(R))=(\tt{1},\tt{011001001100100110},\tt{1},\tt{0110010})$
    is a $\tfrac{7}{3}$-power.  Next we find that $\varphi^2(R)=(\tt{2},z,\tt{2},x)$,
    where $z=\tt{0220102010220102010220}$ and $x=\tt{022010201}$.  We have
    \begin{align*}
        |z|_\tt{2}=8&>7=\tfrac{7}{3}|x|_\tt{2},\\
        |z|_\tt{1}+|z|_\tt{2}=12&> \tfrac{35}{3} = \tfrac{7}{3}\left(|x|_\tt{1}+|x|_\tt{2}\right), \text{ and}\\
        |z|_\tt{0}+|z|_\tt{1}+|z|_\tt{2}=22&>21=\tfrac{7}{3}\left(|x|_\tt{0}+|x|_\tt{1}+|x|_\tt{2}\right),
    \end{align*}
    so $\varphi^2(R)$ is a top-heavy $\tfrac{7}{3}$-power.  
    
    To summarize, we have shown that $\rho(R)$ and $\rho(\varphi(R))$ are $\tfrac{7}{3}$-powers, and that $\varphi^2(R)$ is a top-heavy $\tfrac{7}{3}$-power.  By the \nameref{Proposition:Persistence}, we conclude that $R$ is $\tfrac{7}{3}$-persistent with respect to $\varphi$ and $\rho$ with tipping point $2$.
\end{example}

\section{The structure of \texorpdfstring{$\tfrac{16}{7}$}{16/7}-power-free ternary rich words}
\label{Section:TernaryStructure}

Throughout this section, let $\varphi=\varphi_5$, and let $\pi:\Sigma_5^*\rightarrow \Sigma_3^*$ denote the morphism defined in \autoref{Theorem:TernaryStructure}, i.e., we have
\begin{align*}
    \pi(\tt{0})&=A\\
    \pi(\tt{1})&=AD\\
    \pi(\tt{2})&=AC\\
    \pi(\tt{3})&=ACC\\
    \pi(\tt{4})&=ACCBCC,
\end{align*}    
where $A=\tt{00101101}$, $B=\tt{001}$, $C=\tt{00202202}$, and $D=\tt{002}$.  Finally, let $\hat{\pi}=\pi\circ \varphi^2$. 
Observe that $\hat{\pi}$ is in class $P_\mathrm{ret}$ with marker $p=AADAA\tt{00}$, and that $|\hat{\pi}(\tt{0})|\leq |\hat{\pi}(\tt{1})|\leq \cdots\leq |\hat{\pi}(\tt{4})|$.  The goal of this section is to prove the following statement, from which the structural part of \autoref{Theorem:TernaryStructure} follows easily.

\begin{theorem}\label{Theorem:TernaryStructureInduction}
    Let $\zz$ be a recurrent $\tfrac{16}{7}$-power-free ternary rich word. Then for all $n\geq 0$, up to permutation of the letters, some suffix of $\zz$ has the form $\hat{\pi}(\varphi^n(\xx))$ for some infinite word $\xx\in\Sigma_5^\omega$.
\end{theorem}

We will prove this statement by induction on $n$.  We first prove two lemmas which will make the proof easier to read.  The first lemma establishes the base case, and the second lemma is useful at the inductive step.  We note that although we established a~weak ``outer layer'' of \autoref{Theorem:TernaryStructure} in \autoref{Example:OuterLayerForTernaryRichWords}, we actually need to establish a slightly stronger outer layer in order for our proof of the inductive step of \autoref{Theorem:TernaryStructureInduction} to work.  This is why we use $\hat{\pi}$ in place of $\pi$ in the statement of \autoref{Theorem:TernaryStructureInduction}.  Essentially, we roll two iterations of $\varphi$ into the ``outer layer''.

\begin{lemma}\label{Lemma:TernaryBaseCase}
    Let $\zz$ be a recurrent $\tfrac{16}{7}$-power-free ternary rich word.  Then up to permutation of the letters, some suffix of $\zz$ has the form $\hat{\pi}(\xx)=\pi(\varphi^2(\xx))$ for some infinite word $\xx\in\Sigma_5^\omega$.
\end{lemma}

\begin{proof}
    A backtracking search shows that there are only finitely many $\tfrac{16}{7}$-power-free ternary rich words containing no permutation of the word $P=AADAACCBCCAD$.  (The longest such word has length 1348.)  So by permuting the letters of $\zz$ if necessary, we may assume that $\zz$ contains the factor $P$.  Now by a backtracking search for complete return words to $p=AADAA\tt{00}$ in $\zz$ compatible with the factor $P$, we find that
    \[
        \mathcal{R}_\zz(p)\subseteq \{\hat{\pi}(c)\colon\ c\in\Sigma_5\}.
    \]
    Thus, up to permutation of the letters, some suffix of $\zz$ has the form $\hat{\pi}(\xx)$ for some word $\xx\in\Sigma_5^\omega$.
\end{proof}

\begin{lemma}\label{Lemma:TernaryInductiveStep}
    Let $\xx\in\Sigma_5^\omega$ be a recurrent word such that $\varphi^2(\xx)$ is rich and $\xx$ contains no $\tfrac{16}{7}$-persistent bookended repetition (with respect to $\varphi$ and $\hat{\pi}$).  Then some suffix of $\xx$ has the form $\varphi(\yy)$ for some infinite word $\yy\in\Sigma_5^\omega$.
\end{lemma}

\begin{proof}
    First of all, a backtracking search shows that there are only finitely many words $w$ over $\{\tt{1,2,3,4}\}$ such that $\varphi^2(w)$ is rich and $w$ contains no $\tfrac{16}{7}$-persistent bookended repetition.  (The longest such word has length $3$.)  It follows that $\xx$ contains the letter $\tt{0}$.
    
    It now suffices to show that $\mathcal{R}_\xx(\tt{0})\subseteq \{\tt{01},\tt{02},\tt{03},\tt{04},\tt{044}\}$.  We first show that $\mathcal{R}_\xx(\tt{0})$ is a subset of a slightly larger set.  If $w$ is a complete return word to $\tt{0}$ in $\xx$, then we know by our assumptions about $\xx$ that
    \begin{itemize}
        \item $\varphi^2(w)$ is rich (hence $w$ is also rich as $\varphi$ is in class $P_\mathrm{ret}$), and
        \item $w$ contains no $\tfrac{16}{7}$-persistent bookended repetition with tipping point at most $100$.
    \end{itemize}
    By a depth-first search for complete return words $w$ to $\tt{0}$ over $\Sigma_5$ satisfying these two conditions, we find that $\mathcal{R}_\xx(\tt{0})$ must be a subset of
    \[
        S=\{\tt{01,0121,0131,0141,02,0232,0242,03,0343,04,044}\}.
    \]
    For a more detailed description/illustration of this search, see \autoref{Example:ReturnsTo0}.
    
    Recall \autoref{Definition:Compatible}, and observe that some pairs of return words in $S$  are compatible (e.g., $\tt{01}$ and $\tt{0121}$) while others are not (e.g., $\tt{0121}$ and $\tt{02}$).  By \autoref{Observation:CompatibleSets}, the set $\mathcal{R}_\xx(\tt{0})$ must be compatible, hence it must be contained in some maximal compatible subset of $S$.  Using a computer, we find that there are $15$ maximal compatible subsets of $S$:    
    \begin{align*}
        S_1&=\{\tt{01,0121,0131,0141}\} & S_9&=\{\tt{01,02,03,0242}\}\\
        S_2&=\{\tt{01,0121,03,0141}\} & S_{10}&=\{\tt{01,02,03,0343}\}\\
        S_3&=\{\tt{01,0121,03,0343}\} & S_{11}&=\{\tt{01,0121,0131,04,044}\}\\
        S_4&=\{\tt{01,02,0131,0141}\} & S_{12}&=\{\tt{01,0121,03,04,044}\}\\
        S_5&=\{\tt{01,02,131,0242}\} & S_{13}&=\{\tt{01,02,0131,04,044}\}\\
        S_6&=\{\tt{01,02,0232,0141}\} & S_{14}&=\{\tt{01,02,0232,04,044}\}\\
        S_7&=\{\tt{01,02,03,0141}\} & S_{15}&=\{\tt{01,02,03,04,044}\}\\
        S_8&=\{\tt{01,02,0232,0242}\}
    \end{align*}
    For each $i=1,2,\ldots,15$, the set $S_i$ corresponds to a morphism $\psi_i$ in the obvious manner, e.g., the set $S_1$ corresponds to the morphism $\psi_1:\Sigma_4^*\rightarrow \Sigma_5^*$ defined by 
    \begin{align*}
    \psi_1(\tt{0})&=\tt{01}\\
    \psi_1(\tt{1})&=\tt{0121}\\
    \psi_1(\tt{2})&=\tt{0131}\\
    \psi_1(\tt{3})&=\tt{0141}.
    \end{align*}
    In particular, note that we have $\psi_{15}=\varphi$, and that all of the $\psi_i$'s are in class $P_\mathrm{ret}$ with marker $\tt{0}$.
    
    So for some $i\in\{1,2,\ldots,15\}$, we have  $\mathcal{R}_\xx(\tt{0})\subseteq S_i$, which means that some suffix of $\xx$ has the form $\psi_i(\yy)$ for some $\yy\in\Sigma_5^\omega$.  We claim that we must have $i=15$.  Suppose otherwise that  $\psi_i(\yy)$ is a suffix of $\xx$ for some $i\in\{1,2,\ldots,14\}$.  Then $\varphi^2(\psi_i(\yy))$ is rich and $\yy$ contains no bookended repetition $R$ such that $\psi_i(R)$ is $\tfrac{16}{7}$-persistent.  But by a separate backtracking search for each $i\in\{1,2,\ldots,14\}$, we find that there are only finitely many words $w$ such that
    \begin{enumerate}[label=(\roman*)]
        \item \label{Item:Condition1Ternary} $\varphi^2(\psi_i(w))$ is rich, and
        \item \label{Item:Condition2Ternary} 
        $w$ contains no bookended repetition $R$ such that $\psi_i(R)$ is $\tfrac{16}{7}$-persistent with tipping point at most $100$.
    \end{enumerate}
    For each $i\in\{1,2,\ldots,14\}$, the length of the longest word $w$ satisfying conditions \ref{Item:Condition1Ternary} and \ref{Item:Condition2Ternary} is given in \autoref{Table:PsiBacktracking}.  For a more detailed description/illustration of the search in the case $i=4$, see \autoref{Example:PsiBacktracking}.

    \begin{table}
        \centering
        \begin{tabular}{c|cccccccccccccc}
             $i$ & 1 & 2 & 3 & 4 & 5 & 6 & 7 & 8 & 9 & 10 & 11 & 12 & 13 & 14\\\hline
             
             \makecell{Length of a longest word sat- \\
             isfying conditions \ref{Item:Condition1Ternary} and \ref{Item:Condition2Ternary}} &
             8 & 4 & 9 & 3 & 7 & 7 & 8 & 11 & 8 & 29 & 9 & 15 & 34 & 34
        \end{tabular}
        \caption{Results of backtracking searches in the proof of \autoref{Lemma:TernaryInductiveStep}.}
        \label{Table:PsiBacktracking}
    \end{table}

    We conclude that some suffix of $\yy$ has the form $\psi_{15}(\xx)=\varphi(\xx)$, i.e., that $\mathcal{R}_\xx(\tt{0})\subseteq \{\tt{01,02,03,04,044}\}$.
\end{proof}

\begin{proof}[Proof of \autoref{Theorem:TernaryStructureInduction}]
    We proceed by induction on $n$.  By \autoref{Lemma:TernaryBaseCase}, the statement holds when $n=0$.  Now suppose for some integer $k\geq 0$ that some suffix of $\zz$ has the form $\hat{\pi}(\varphi^k(\xx))$, where $\xx\in\Sigma_5^\omega$.  Since $\hat{\pi}(\varphi^k(\xx))=\pi(\varphi^{k+2}(\xx))$ is rich and both $\pi$ and $\varphi$ belong to class $P_\mathrm{ret}$, we know that $\varphi^2(\xx)$ is rich.  Further, since $\hat{\pi}(\varphi^k(\xx))$ is $\tfrac{16}{7}$-power-free, we know that $\xx$ contains no $\tfrac{16}{7}$-persistent bookended repetition (with respect to $\varphi$ and $\hat{\pi}$).  Thus, by \autoref{Lemma:TernaryInductiveStep}, some suffix of $\xx$ has the form $\varphi(\yy)$ for some word $\yy\in\Sigma_5^\omega$.  It follows that $\zz$ has a suffix of the form $\hat{\pi}(\varphi^{k+1}(\yy))$, and this completes the inductive step. 
\end{proof}

\begin{corollary}  The language of $\pi(\uu_{5})$ is contained in the language of every recurrent $\tfrac{16}{7}$-power-free ternary rich word $\zz$. Hence,   $RT(\mathcal{C}_{3}) = E\bigl({\pi}(\uu_{5})\bigr)$.
\end{corollary}

\begin{proof} By \autoref{Theorem:TernaryStructureInduction}, for each $n \in \N, $ the language of $\zz$ contains every factor of $\hat{\pi}(\varphi^{n}(\xx))= \pi(\varphi^{n+2}(\xx))$ for some infinite word $\xx\in\Sigma_5^\omega$.  Since  the letter $\tt 0$ appears in $\varphi(\xx)$, the factor ${\pi}\bigl(\varphi^n(\tt 0)\bigl)$ appears in $\zz$ for every $n\in \N$. The corollary follows from the facts  that  ${\pi}(\uu_5)$ is rich, $E\bigl({\pi}(\uu_{5})\bigr)<\tfrac{16}{7}$ and $\uu_5= \varphi^\omega(\tt 0)$.
\end{proof}

\section{The structure of \texorpdfstring{$\tfrac{53}{25}$}{53/25}-power-free quaternary rich words}
\label{Section:QuaternaryStructure}

Throughout this section, let $\varphi=\varphi_7$, and let $\pi_1,\pi_2,\ldots,\pi_7:\Sigma_7^*\rightarrow \Sigma_4^*$ denote the morphisms defined in \autoref{Table:PiDefinition}.
For each $i\in\{1,2,\ldots,7\}$, let $\hat{\pi}_i=\pi_i\circ \varphi^4$. For each $i\in\{1,2,\ldots,7\}$, since $\pi_i$ and $\varphi$ are both in class $P_{\mathrm{ret}}$, we see that $\hat{\pi}_i$ is also in class $P_\mathrm{ret}$, and it is straightforward to check that $|\hat{\pi}_i(\tt{0})|\leq |\hat{\pi}_i(\tt{1})|\leq \cdots\leq |\hat{\pi}_i(\tt{6})|$.  We will prove the following statement, from which the structural part of \autoref{Theorem:QuaternaryStructure} follows easily.

\begin{theorem}\label{Theorem:QuaternaryStructureInduction}
    Let $\zz$ be a recurrent $\tfrac{53}{25}$-power-free quaternary rich word. Then for all $n\geq 0$, up to permutation of the letters, some suffix of $\zz$ has the form $\hat{\pi}_i(\varphi^n(\xx))$ for some $i\in\{1,2,\ldots,7\}$ and some infinite word $\xx\in\Sigma_7^\omega$.
\end{theorem}

Our overall approach is similar to the one we used to prove \autoref{Theorem:TernaryStructureInduction} in \autoref{Section:TernaryStructure}.  It is complicated by the fact that we must deal with seven different projections $\hat{\pi}_1,\hat{\pi}_2,\ldots,\hat{\pi}_7$, and by the fact that the backtracking searches involved are much more time-consuming.

We first establish the ``outer layer'', which will serve as the base case of our inductive proof of \autoref{Theorem:QuaternaryStructureInduction}.  We use the following ``lexicographic trick'' to speed up some of the backtracking searches in the proof.  Suppose that $L$ is a factorial language, and that we wish to show that there are only finitely many words in $L$ with prefix $p$.  Our search traverses the tree of all words in $L$ with prefix $p$ in lexicographic order, backtracking whenever we encounter a word not in $L$.  Now suppose that at some point during this search, we encounter a word $u$ that contains a suffix $v$ such that $v$ has prefix $p$ and $v<u$ in the lexicographic order.  Since we are searching through the tree in lexicographic order, we already know that there are only finitely many words in $L$ with prefix $v$ (as we have already made it past $v$ to $u$).  Hence, since $v$ is a factor of $u$, we know that there are only finitely many words in $L$ with prefix $u$.  Thus, when we encounter such a word $u$ in the tree, we may backtrack immediately.  Essentially, this allows us to avoid going down some ``dead ends'' in the tree.  When we use this lexicographic trick, we do not necessarily find a longest word in $L$ with prefix $p$; instead we report the length of a longest word encountered in the search.

\begin{lemma}\label{Lemma:QuaternaryBaseCase}
    Let $\zz$ be an infinite $\tfrac{53}{25}$-power-free quaternary rich word.  Then up to permutation of the letters, some suffix of $\zz$ has the form $\hat{\pi}_i(\xx)=\pi_i(\varphi^4(\xx))$ for some $i\in\{1,2,\ldots,7\}$ and some infinite word $\xx\in\Sigma_7^\omega$.
\end{lemma}

\begin{proof}
We first claim that $\zz$ contains some permutation of $\hat{\pi}_i(\tt{5})$ for some $i\in\{1,2,\ldots,7\}$.  In principle, this claim could be established with a single backtracking search, but this approach proved to be too time-consuming, so we broke up the search as follows.
 
A backtracking search shows that a longest $\tfrac{53}{25}$-power-free rich word over $\Sigma_4$ containing no factor from the set $\{\tt{00},\tt{11},\tt{22},\tt{33}\}$ has length $22$.  So by permuting the letters if necessary, we may assume that $\zz$ has prefix $\tt{001}$.  Next, a backtracking search using the ``lexicographic trick'' described above shows that there are only finitely many $\tfrac{53}{25}$-power-free rich words over $\Sigma_4$ with prefix $\tt{001}$ containing no permutation of the following seven words:
\begin{align*}
 w_1&=\tt{0010011001002}\\
 w_2&=\tt{010010110100102}\\
 w_3&=\tt{01011010010110102}\\
 w_4&=\tt{011001102}\\
 w_5&=\tt{0101101022}\\
 w_6&=\tt{00101101002}\\
 w_7&=\tt{0101100110102}
\end{align*}
A longest word encountered in the search has length $1951$.  So we may assume that $\zz$ contains $w_i$ as a factor for some $i\in\{1,\ldots,7\}$.  

Now we aim to show that for each $i\in\{1,2,\ldots,7\}$, if $\zz$ has prefix $w_i$, then $\zz$ contains some permutation of $\hat{\pi}_i(\tt{5})$.  We use seven separate backtracking searches with the ``lexicographic trick'' to establish that for each $i\in\{1,2,\ldots,7\}$, there are only finitely many $\tfrac{53}{25}$-power-free rich words over $\Sigma_4$ with prefix $w_i$ that contain no permutation of $\hat{\pi}_i(\tt{5})$.  The length of a longest word encountered in each of the seven searches is reported in~\autoref{Table:Backtracking}. 
So by permuting the letters and taking a suffix if necessary, we may assume that $\mathbf{z}$ has prefix $\hat{\pi}_i(\tt{5})$ for some $i\in\{1,2,\ldots,7\}$.

\begin{table}
    \centering
    \begin{tabular}{c|c}
         $i$ & \makecell{Length of a longest word encountered}\\\hline
         1 & 26290\\
         2 & 6883\\
         3 & 7959\\
         4 & 74978\\
         5 & 24253\\
         6 & 36772\\
         7 & 45458
    \end{tabular}
    \caption{Results of our backtracking searches for $\tfrac{53}{25}$-power-free rich words over $\Sigma_4$ with prefix $w_i$ containing no permutation of $\hat{\pi}_i(\tt{5})$.}
    \label{Table:Backtracking}
\end{table}
    
    Fix $i\in\{1,2,\ldots,7\}$ and suppose that $\zz$ contains $\hat{\pi}_i(\tt{5})$.  Recall that $\hat{\pi}_i$ is in class $P_\mathrm{ret}$, and let $p_i$ denote the corresponding marker.  A backtracking search for complete return words to $p_i$ in $\zz$ compatible with $\hat{\pi}_i(\tt{5})$ tells us that
    \[
        \mathcal{R}_\zz(p_i)\subseteq \{\hat{\pi}_i(c)\colon\ c\in\Sigma_7\}.
    \]
    Thus, some suffix of $\zz$ has the form $\hat{\pi}_i(\xx)$ for some word $\xx\in\Sigma_7^\omega$.    
\end{proof}

\begin{lemma}\label{Lemma:QuaternaryInductiveStep}
    For each $i\in\{1,2,\ldots,7\}$, if $\xx\in\Sigma_7^\omega$ is a recurrent word such that $\varphi^4(\xx)$ is rich and $\xx$ contains no $\tfrac{53}{25}$-persistent bookended repetition with respect to $\varphi$ and $\hat{\pi}_i$, then some suffix of $\xx$ has the form $\varphi(\yy)$ for some infinite word $\yy\in\Sigma_7^\omega$.
\end{lemma}

\begin{proof}
    We have seven separate cases corresponding to $i=1,2,\ldots, 7$, but they all proceed in a similar manner.  Fix $i\in\{1,2,\ldots,7\}$, and let $\xx\in\Sigma_7^\omega$ be a recurrent word such that $\varphi^4(\xx)$ is rich and $\xx$ contains no $\tfrac{53}{25}$-persistent bookended repetition with respect to $\varphi$ and $\hat{\pi}_i$.  First of all, a backtracking search shows that there are only finitely many words $w$ over $\{\tt{1,2,3,4,5,6}\}$ such that $\varphi^4(w)$ is rich and $w$ contains no $\tfrac{53}{25}$-persistent bookended repetition with respect to $\varphi$ and $\hat{\pi}_i$ with tipping point at most $100$.  It follows that $\xx$ contains the letter $\tt{0}$.
    
    It now suffices to show that $\mathcal{R}_\xx(\tt{0})\subseteq \{\varphi(c)\colon\ c\in\Sigma_7\}=\{\tt{01},\tt{02},\tt{03},\tt{04},\tt{05},\tt{06},\tt{066}\}$.  We first show that $\mathcal{R}_\xx(\tt{0})$ is a subset of a slightly larger set.  If $w$ is a complete return word to $\tt{0}$ in $\xx$, then we know by our assumptions about $\xx$ that
    \begin{itemize}
        \item $\varphi^4(w)$ is rich (hence $w$ is also rich as $\varphi$ is in class $P_\mathrm{ret}$), and
        \item $w$ contains no $\tfrac{53}{25}$-persistent bookended repetition with respect to $\varphi$ and $\hat{\pi}_i$ with tipping point at most $100$.
    \end{itemize}
    By a depth-first search for complete return words $w$ to $\tt{0}$ over $\Sigma_5$ satisfying these two conditions, we find that $\mathcal{R}_\xx(\tt{0})$ must be a subset of
    \begin{align*}
        S&=\{\tt{01,0121,0131,0141,0151,0161,02,0232,0242,0252,0262},\\
        &\ \ \ \ \ \, \tt{03,0343,0353,0363,04,0454,0464,05,0565,06,066}\}.
    \end{align*}
    To be completely clear, we find this same set $S$ for each $i\in\{1,2,\ldots,7\}$.\footnote{This was slightly surprising to us, but after all, the maps $\hat{\pi}_1,\hat{\pi}_2,\ldots,\hat{\pi}_7$ do have similar structure.}
    
    By \autoref{Observation:CompatibleSets}, the set $\mathcal{R}_\xx(\tt{0})$ must be compatible, hence it must be contained in some maximal compatible subset of $S$.  Using a computer, we find that there are $203$ maximal compatible subsets of $S$ (each of them of size at most seven), which we denote by $S_1,S_2,\ldots,S_{203}$.  We do not write all of these sets explicitly here, but the interested reader can see them by running our code.  We note that $S_{203}=\{\tt{01},\tt{02},\tt{03},\tt{04},\tt{05},\tt{06},\tt{066}\}=\{\varphi(c)\colon\ c\in\Sigma_7\}$.
    
    For each $j=1,2,\ldots,203$, the set $S_j$ corresponds to a morphism $\psi_j$ in the obvious manner.  In particular, we have $\psi_{203}=\varphi$, and we note that all of the $\psi_j$'s are in class $P_\mathrm{ret}$ with marker $\tt{0}$.
    So for some $j\in\{1,2,\ldots,203\}$, we have  $\mathcal{R}_\xx(\tt{0})\subseteq S_j$, which means that some suffix of $\xx$ has the form $\psi_j(\yy)$ for some $\yy\in\Sigma_7^\omega$.  We claim that we must have $j=203$.  Suppose otherwise that  $\psi_j(\yy)$ is a suffix of $\xx$ for some $j\in\{1,2,\ldots,202\}$.  Then $\varphi^4(\psi_j(\yy))$ is rich and $\yy$ contains no bookended repetition $R$ such that $\psi_j(R)$ is $\tfrac{53}{25}$-persistent with respect to $\varphi$ and $\hat{\pi}_i$.  But by a separate backtracking search for each $j\in\{1,2,\ldots,202\}$, we find that there are only finitely many words $w$ such that
    \begin{enumerate}[label=(\roman*)]
        \item \label{Item:Condition1Quaternary} $\varphi^4(\psi_j(w))$ is rich, and
        \item \label{Item:Condition2Quaternary} 
        $w$ contains no bookended repetition $R$ such that $\psi_j(R)$ is $\tfrac{53}{25}$-persistent with respect to $\varphi$ and $\hat{\pi}_i$ with tipping point at most 100.
    \end{enumerate}
    The interested reader can see more details concerning the results of these backtracking searches by running our code.

    We conclude that some suffix of $\yy$ has the form $\psi_{203}(\xx)=\varphi(\xx)$, i.e., that $\mathcal{R}_\xx(\tt{0})\subseteq \{\tt{01,02,03,04,05,06,066}\}$, as desired.
\end{proof}

\begin{proof}[Proof of \autoref{Theorem:QuaternaryStructureInduction}]
    We proceed by induction on $n$.  By \autoref{Lemma:QuaternaryBaseCase}, the statement holds when $n=0$.  Now suppose for some integer $k\geq 0$ that some suffix of $\zz$ has the form $\hat{\pi}_i(\varphi^k(\xx))$ for some $i\in\{1,2,\ldots,7\}$ and some word $\xx\in\Sigma_7^\omega$.  Since $\hat{\pi}_i(\varphi^k(\xx))=\pi_i(\varphi^{k+4}(\xx))$ is rich and both $\pi_i$ and $\varphi$ belong to class $P_\mathrm{ret}$, we know that $\varphi^4(\xx)$ is rich.  Further, since $\hat{\pi}_i(\varphi^k(\xx))$ is $\tfrac{53}{25}$-power-free, we know that $\xx$ contains no $\tfrac{53}{25}$-persistent bookended repetition with respect to $\varphi$ and $\hat{\pi}_i$.  Thus, by \autoref{Lemma:QuaternaryInductiveStep}, some suffix of $\xx$ has the form $\varphi(\yy)$ for some word $\yy\in\Sigma_7^\omega$.  It follows that $\zz$ has a suffix of the form $\hat{\pi}_i(\varphi^{k+1}(\yy))$, and this completes the inductive step. 
\end{proof}

\section{Richness and Critical Exponent of the Words \texorpdfstring{$\pi_i(\uu_7)$}{pi(u)}}\label{Section:Pi1}

In this section, 
we study the word $\zz = \pi(\uu)$, where $\uu=\uu_{7}$ is  the fixed point of the morphism $\varphi=\varphi_7$ over the alphabet 
$\Sigma_7$ and $\pi=\pi_1$ is the projection from $\Sigma_7^*$ to $\Sigma_4^*$ defined in \autoref{Section:Introduction}.  We will show that $\zz$ is rich (and moreover $G$-rich for a group $G$ generated by two involutive antimorphisms), and that $\zz$ has critical exponent equal to $E^*(\uu_7)$.  This completes the proof of \autoref{con:opet} for $d=4$ and the proof of \autoref{Theorem:QuaternaryStructure} for $d=4$ in the case $i=1$.  The remaining cases $i=2,3,\ldots,7$ can be established by the same method, and the details are only slightly different in each case, so we omit the proofs of the other cases.

For convenience, we recall that  
\begin{equation}\label{eq:morphismFi7}\varphi: \ \begin{array}{l}
    \tt{0} \mapsto \tt{01}\\
    \tt{1} \mapsto \tt{02}\\
    \tt{2} \mapsto \tt{03}\\
    \tt{3} \mapsto \tt{04}\\
    \tt{4} \mapsto \tt{05}\\
    \tt{5} \mapsto \tt{06}\\
    \tt{6} \mapsto \tt{066} 
\end{array}
\ \ \qquad \text{and}\qquad \ \  
\pi: \ \begin{array}{l}
    \tt{0} \mapsto \tt{0}\\
    \tt{1} \mapsto \tt{01}\\
    \tt{2} \mapsto \tt{011}\\
    \tt{3} \mapsto \tt{02}\\
    \tt{4} \mapsto \tt{022}\\
    \tt{5} \mapsto \tt{022322}\\
    \tt{6} \mapsto \tt{0223223322322}.
\end{array}
\end{equation}
Here are prefixes of $\uu$ and $\zz$:
$$\uu=\tt{0102010301020104010201030102010501020103010201040102010301020106}\cdots$$
$$\zz=\tt{0010011001002001001100100220010011001002001001100100223220010011}\cdots$$
First we will prove in \autoref{subsec:symmetries} that  $\zz$ is rich and then we will study bispecial factors in $\uu$ and $\zz$  with the aim of determining its critical exponent as expressed in \autoref{thm:main} below.   

\begin{theorem}\label{thm:main} Let $\zz = \pi_1(\uu_7)$. Then 
$$E({\bf{z}})=E^*(\uu_7)=1+\frac{1}{3-\Lambda}\doteq 2.1197196907\,,$$ 
where $\Lambda$ is the root of $x^4-2x^3-1$ of largest modulus.
\end{theorem}
The main tool that we use to find the critical exponent of $\zz$ is the following theorem. 

\begin{theorem}[{\cite[Theorem 3]{DolceDP2023}}]\label{thm:FormulaForE}
Let $\xx$ be a uniformly recurrent aperiodic word.
Then
$$
E(\xx) = 1 + \sup\left\{\tfrac{|w|}{|r|}\, : \, \text{$w$ is a bispecial factor of $\xx$ and $r$ is a return word to $ w$ in $\xx$ }\right\}\,.
$$
\end{theorem}

Let us mention that using exactly the same procedure, one can show that $\pi_i(\uu_7)$ is rich and has the same critical exponent as $\zz$ for all $i \in \{2,\ldots, 7\}$. Further, we note that the languages of $\pi_1(\uu_7), \ldots, \pi_7(\uu_7)$ are mutually distinct.  To see this, we observe that return words to letters in each of these words do not coincide even after any letter permutation.

\subsection{Richness of \texorpdfstring{$\zz=\pi_1(\uu_7)$}{z=pi1(u7)}}\label{subsec:symmetries}

We begin by studying the symmetries that the infinite word $\zz$ possesses. These symmetries will help us to establish the richness of $\zz$, and will also make it easier to investigate the critical exponent of $\zz$.

\begin{proposition}\label{prop:symmetry}
The language of $\zz$ is closed under $R$ and $S$, where $R$ is the reversal mapping and $S$ is the antimorphism swapping $\tt{0}\leftrightarrow \tt{2}$ and $\tt{1}\leftrightarrow \tt{3}$. 
\end{proposition}
\begin{proof}
As $\uu$  is fixed by a primitive morphism, its language is uniformly recurrent. Hence the language of $\zz= \pi(\uu)$ is uniformly recurrent. As $\uu$ is rich~\cite{DvorakovaKloudaPelantova2025},  $\uu$ contains  infinitely many palindromes and  is closed under reversal. Since $\pi$ is of class $P_{\mathrm{ret}}$, the language of  $\pi(\uu)$ is closed under reversal as well.

Closedness under $S$ will be clear if we show that $\zz$ contains $S$-images of factors of unbounded length. It suffices to show that for all $n \in \mathbb N$
\begin{equation}\label{eq:obrazByS}S(\pi(\varphi^{n}(\tt{6}))) =  \pi(\varphi^{n+3}(\tt{0}))\,.\end{equation}
Let us prove it by induction. 
For $n=0$, we have 
\[
S(\pi(\tt{6}))=S(\tt{0223223322322})=\tt{0010011001002}
\]
and
\[
\pi(\varphi^3(\tt{0}))=\pi(\tt{01020103})=\tt{0010011001002}\,,
\]
giving the desired equality.  For $n\in \{1,2\}$, we invite the reader to check the equality.

For $n\geq 3$, we obtain
\begin{align*}
    S(\pi(\varphi^{n+1}(\tt{6})))
&=S(\pi(\varphi^{n}(\tt{066})))\\
&=S(\pi(\varphi^n(\tt{6})))S(\pi(\varphi^n(\tt{6})))S(\pi(\varphi^n(\tt{0})))\\
&= \pi(\varphi^{n+3}(\tt{0})) \pi(\varphi^{n+3}(\tt{0})) \pi(\varphi^{n-3}(\tt{6}))\\
&= \pi(\varphi^{n+3}(\tt{0})) \pi(\varphi^{n+2}(\tt{0})\varphi^{n+1}(\tt{0}) \cdots \varphi^{n-3}(\tt{0})\varphi^{n-3}(\tt{6}))\pi(\varphi^{n-3}(\tt{6}))\,.
\end{align*}
where on the second line, we used the induction assumption that $S(\pi(\varphi^{n}(\tt{6}))) =  \pi(\varphi^{n+3}(\tt{0}))$ and $\pi(\varphi^{n-3}(\tt{6})) =  S(\pi(\varphi^{n}(\tt{0})))$, and besides that we repeatedly used the definition of $\varphi$ and the fact that $S$ is an antimorphism.
Similarly, using the definition of $\varphi$, we derive
\begin{align*}
\pi(\varphi^{n+4}(\tt{0}))&= \pi(\varphi^{n+3}(\tt{0})) \pi(\varphi^{n+3}(\tt{1}))\\
&=\pi(\varphi^{n+3}(\tt{0}))\pi(\varphi^{n+2}(\tt{0})\varphi^{n+1}(\tt{0})\cdots \varphi^{n-3}(\tt{0})\varphi^{n-3}(\tt{6})\varphi^{n-3}(\tt{6}))\,. \qedhere
\end{align*}
\end{proof}

Since the projection $\pi$ is of class $P_{\mathrm{ret}}$,  the proof of the richness of $\zz$ can be simplified by using the following statement based on bilateral orders of bispecial factors. 

\begin{proposition}\label{prop:richnessPret}
Let $\eta: \Sigma_d^*\to \Sigma_k^*$ be a morphism of class $P_{\mathrm{ret}}$ from \autoref{eq:Pret}. If every bispecial factor $w$ in $\eta(\uu_d)$, which does not contain the palindromic marker $p$, satisfies
\begin{itemize}
\item ${\mathrm b}(w)=0$ if $w$ is not a palindrome, and
\item ${\mathrm b}(w)=\#{\mathrm{Pext}}(w)-1$ if $w$ is a palindrome,
\end{itemize}
then $\eta(\uu_d)$ is rich.
\end{proposition}
\begin{proof}
This statement was proven as Proposition 19 in \cite{DvPe_LimRT2025} for a more restricted definition of class $P_{\mathrm{ret}}$, where all letter images take on the form $\eta(i)=pq_i$ with $p, q_i$ palindromes. However, the same proof applies here, too.
\end{proof}

\begin{theorem}\label{thm:rich}
The infinite word $\zz=\pi_1(\uu_7)$ is rich. 
\end{theorem}
\begin{proof}
Using \autoref{prop:richnessPret}, it suffices to find bilateral orders of bispecial factors of $\zz$ that do not contain $\tt 0$. 
Moreover, using \autoref{prop:symmetry}, we can make use of closedness of the language of $\zz$ under $S$.
Obviously, if $w$ is a bispecial factor in $\zz$, then ${\mathrm b}(w)={\mathrm b{(S(w))}}$ and $\#\mathrm{Pext}(w)=\#\mathrm{Pext}(S(w))$. 
Therefore, it even suffices to find bilateral orders of bispecial factors of $\zz$ that contain neither $\tt{0}$ nor $\tt{2}$.
The unique such bispecial factors are $\tt 1$ and $\tt 3$. 
The bilateral order of $\tt 1$ satisfies ${\rm b}({\tt 1})=0=\#{\rm Pext}({\tt 1})-1$ since
$${\rm Bext}({\tt 1})=\{\tt{010, 011, 110}\}, \ {\rm Pext}({\tt 1})=\{\tt{010}\}\,.$$
Since $S(\tt{1}) = \tt 3$, the bilateral order of $\tt 3$ satisfies ${\rm b}({\tt 3})=0=\#{\rm Pext}({\tt 3})-1$. 
Hence the richness of $\zz$ is a consequence of \autoref{prop:richnessPret}.
    \end{proof}

We will now show that the infinite word $\zz$ is rich in ``generalized palindromes''.  This is necessary for the proof of \autoref{Theorem:QuaternaryStructure} and it supports \autoref{con:structure} saying that rich words with the minimum critical exponent are rich in this more general sense.

First we must describe this generalized concept of richness more precisely.  For an infinite word $\xx$, we let
\begin{itemize}
    \item $\mathcal{C}(n)$ denote the number of factors of $\xx$ of length $n$,
    \item $\mathcal{P}(n)$ denote the number of factors of $\xx$ of length $n$ that are palindromes, and
    \item $\mathcal{P}_S(n)$ denote the number of factors of $\xx$ of length $n$ that are $S$-palindromes (i.e., words $w$ satisfying $S(w)=w$).
\end{itemize}
The functions $\mathcal{C}$, $\mathcal{P}$, and $\mathcal{P}_S$ are called the \emph{factor complexity}, \emph{palindromic complexity}, and \emph{$S$-palindromic complexity} of $\xx$, respectively.

By~\cite{BaMaPe2007}, if the language of the infinite word $\xx$ is closed under reversal, then
\begin{equation}\label{eq:richness}
\mathcal{P}(n)+ \mathcal{P}(n+1)\leq  \mathcal{C}(n+1) - \mathcal{C}(n) +2\qquad \text{for every $n \in \N $}.
\end{equation} 
Moreover, by~\cite{BuLuGlZa2009}, the word $\xx$ is rich if and only if equality holds in \eqref{eq:richness} for all $n\in\N$.
By~\cite{PStarosta2013}, if the language of $\xx$ is closed under both $R$ and $S$, then 
\begin{equation}\label{eq:Grichness}\mathcal{P}(n)+ \mathcal{P}(n+1) +\mathcal{P}_S(n)+ \mathcal{P}_S(n+1)   \leq \mathcal{C}(n+1) - \mathcal{C}(n) +4\qquad \text{for every $n \in \N $}.
\end{equation}
The characterization of rich words in terms of \eqref{eq:richness} described above inspires the following definition of $G$-richness for the group $G$ generated by $R$ and $S$: If the language of the infinite word $\xx$ is closed under both $R$ and $S$, then we say that $\xx$ is \emph{$G$-rich} if equality holds in \eqref{eq:Grichness} for all $n\in\N, \ n \geq 1$.

\begin{theorem}\label{thm:Grich} 
The word $\zz=\pi_1(\uu_7)$ is $G$-rich, where $G$ is the group generated by the antimorphisms $R$ and $S$.
\end{theorem}
 
\begin{proof} By \autoref{thm:rich},  $\zz$ is rich and thus $\mathcal{P}(n)+ \mathcal{P}(n+1)=  \mathcal{C}(n+1) - \mathcal{C}(n) +2$ for every $n \in \N$. To demonstrate equality in $\eqref{eq:Grichness}$, we have to show that 
$\mathcal{P}_S(n)+ \mathcal{P}_S(n+1)\geq 2$ for every $n \in \N, \ n\geq 1$. Obviously, any $S$-palindrome has an even length. It means that we have to show that for every $n=2k$ with $k \in \N, \ k\geq 1,$ there exist two $S$-palindromes of length $n$.

Due to \eqref{eq:obrazByS}, the bidirectional infinite word $${\bf y} = \cdots y_{-3}y_{-2}y_{-1}\mid y_1y_2y_3\cdots = \pi\bigr(\varphi^\omega(\tt{6})\bigl)\mid \pi\bigl(\varphi^\omega(\tt{0})\bigr)\,,$$ where $\cdots y_{-3}y_{-2}y_{-1}$ is the left-infinite word having $\pi\bigr(\varphi^n(\tt{6})\bigl)$ as suffix for each $n$, has the property that  $S(\bf y) = \bf y$, i.e.,  the factors   $w=y_{-k}\cdots y_{-1} y_1\cdots y_k$  and $R(w)$ are $S$-palindromes for all $k\in\N$. Since both $w$ and $R(w)$ belong to the language of $\zz$, this completes the proof.  
\end{proof}

\subsection{Bispecial triplets in \texorpdfstring{$\uu ={\bf u}_7$}{u=u7}} \quad 

A~{\em bispecial triplet} in $\uu$ is a~triplet $(\{a,b\}, w, \{c,d\})$, where the pair $\{a, b\}$ is an unordered pair of left extensions of $w$ with $a\not =b$, the pair $\{c,d\}$ is an unordered pair of right extensions of $w$ with  $c\not =d$, and either both $awc$ and $bwd$ or both $awd$ and $bwc$ are factors of $\uu$.  

Clearly, each bispecial factor in $\uu$ is the middle element of some bispecial triplet.   We say that the bispecial factor is \textit{associated} with the bispecial triplet. 
The set of all bispecial triplets in $\uu$ is denoted by $\mathcal{T}$.

Using the method of Klouda~\cite{K2012}, we have shown that every bispecial factor in $\uu$ (or its mirror image) is associated with a~bispecial triplet that can be obtained from a small set of initial bispecial triplets by repeatedly applying the following function $f$.

\begin{definition}\label{def:Zobrazeni_f} 
For $(\{a,\tt{6}\},w,\{b,\tt{6}\})\in  \mathcal{T}$, we set $f(a,w,b) = (\{a',\tt{6}\},w',\{b',\tt{6}\})$, where 
$ a'=(a+1)\bmod 6$,  $b' = (b+1) \bmod 6$, and 
\[
w' =\begin{cases} \phantom{\tt{6}}\varphi(w)\tt 0, & \text{ if } a<5 \text{ and } b<5\,;\\
  \phantom{\tt{6}}\varphi(w)\tt{06}, & \text{ if } a<5 \text{ and } b=5\,;\\

  \tt{6}\varphi(w)\tt{0}, & \text{ if } a=5 \text{ and } b<5\,;\\

  \tt{6}\varphi(w)\tt{06}, & \text{ if } a=5 \text{ and } b=5\,.
  \end{cases}
\]
\end{definition}

Explicitly, with the initial bispecial triplets
\begin{itemize}
\item $T_0=(\{\tt{0},\tt{6}\}, \varepsilon, \{\tt{0},\tt{6}\})$,
\item $T_6=(\{\tt{0},\tt{6}\}, \tt{6}, \{\tt{0},\tt{6}\})$, and
\item $T_k=(\{\tt{0},\tt{6}\}, \varepsilon, \{k,\tt{6}\})$ for $k \in \{\tt{1},\tt{2},\tt{3},\tt{4},\tt{5}\}$,
\end{itemize}
we have that for every bispecial factor $w$ in $\uu$, either the word $w$ or its mirror image $R(w)$ is associated with $f^n(T_i)$ for some $n\geq 0$ and some $i\in \{0,1,\dots, 6\}$.   We note that the bispecial factors associated with $f^n(T_k)$ for $k \not \in \{0, 6\}$ and $n \geq 6-k$ are nonpalindromic, while all other bispecial factors are palindromic.  

\begin{example}\label{ex:triplets} Let us find the first few bispecial triplets obtained from the initial bispecial triplets $T_1$ and $T_2$ by repeatedly applying the function $f$ from  \autoref{def:Zobrazeni_f}. Note that the bispecial factors associated with these initial triplets may coincide, nevertheless the triplets are different because they differ in their right extensions. For  brevity, we use the notation $F_0=\varepsilon$, $F_1 = \tt{0}$, $F_2= \varphi(\tt{0})\tt{0}$,  $F_3=  \varphi^2(\tt{0})\varphi(\tt{0})\tt{0}$, and  $F_4=  \varphi^3(\tt{0})  \varphi^2(\tt{0})\varphi(\tt{0})\tt{0}$.  Note that all of the $F_i$ are palindromes.
\[
\begin{array}{ccl}
  f(T_1)   & = & (\{\tt{1},\tt{6}\}, F_1, \{\tt{2},\tt{6}\})\\
   f^2(T_1)  & = &(\{\tt{2},\tt{6}\}, F_2, \{\tt{3},\tt{6}\}) \\
   f^3(T_1)  & = &(\{\tt{3},\tt{6}\}, F_3, \{\tt{4},\tt{6}\}) \\
      f^4(T_1)  & =& (\{\tt{4},\tt{6}\}, F_4,  \{\tt{5},\tt{6}\}) 
\end{array} \qquad \text{and} \qquad \begin{array}{ccl}
  f(T_2)   & = & (\{\tt{1},\tt{6}\}, F_1, \{\tt{3},\tt{6}\})\\
   f^2(T_2)  & = &(\{\tt{2},\tt{6}\}, F_2, \{\tt{4},\tt{6}\}) \\
   f^3(T_2)  & = &(\{\tt{3},\tt{6}\}, F_3, \{\tt{5},\tt{6}\}) \\
      f^4(T_2)  & =& (\{\tt{4},\tt{6}\}, F_4\tt{6},  \{\tt{0},\tt{6}\}) 
\end{array} 
\]
Finally, we note that  $F_0$, $F_1$, $F_2$, $F_3$, and $F_4$ are the only bispecial factors of $\uu$ with more than two left or right extensions. 
\end{example}

\subsection{Bispecial factors and their shortest return words in \texorpdfstring{$\zz$}{z}}\label{sec:BS}\quad

When searching for all bispecial factors and their shortest return words in $\zz$, we will make use of the fact that the language of $\zz$ is closed under both $R$ and $S$.
\begin{observation}\label{obs:BSsymmetry}
Let $w$ be a bispecial factor in $\zz$ and $r$ its shortest complete return word. Then $R(w)$ and $S(w)$ are bispecial factors in $\zz$ with the shortest complete return words equal to $R(r)$ and  $S(r)$, respectively. 
\end{observation}

On the one hand, bispecial factors in $\uu$ give rise to bispecial factors in $\zz$, as described in the observation below.

\begin{observation}\label{obs:poProjekci}
If $(\{a,b\}, w, \{c,d\})$ is a bispecial triplet in $\uu$, then 
\[
v=\lcs\{\pi(a), \pi(b)\}\,\pi(w)\lcp\{\pi(c),\pi(d)\}
\]
is a bispecial factor in $\zz$ and contains $\tt{0}$ as a factor.
\end{observation}
We will call $v$ the \emph{$\pi$-offspring} of the bispecial triplet $(\{a,b\}, w, \{c,d\})$.

\begin{remark}\label{rem:LcsLcp} The form of $\pi$ implies that for all $i,j$ with $0\leq i<j\leq 6$, we have
\[
\lcs\{\pi(i), \pi(j)\} = \begin{cases}
\lcs\{\pi(i), \pi(\tt{6})\}, & \text{ if $\{i,j\}\neq \{1,2\}$;}\\
\tt{1}, & \text{ if $\{i,j\}=\{1,2\}$;}
\end{cases}
\]
and
\[
\lcp\{\pi(i), \pi(j)\} = 
\tt{0}\lcs\{\pi(i), \pi({j})\}.
\]
\end{remark}

All bispecial factors in $\uu$, except for $F_0$, $F_1$, $F_2$, $F_3$, and $F_4$, have only two left extensions and two right extensions, and have $\tt{6}$ as both a left extension and a right extension. These extensions are included in the bispecial triplets generated from the initial triplets  $$T_0=(\{\tt{0},\tt{6}\}, \varepsilon, \{\tt{0},\tt{6}\}), \quad T_6=(\{\tt{0},\tt{6}\}, \tt{6}, \{\tt{0},\tt{6}\}), \quad T_k=(\{\tt{0},\tt{6}\}, \varepsilon, \{k,\tt{6}\})\quad \text{for} \ 1\leq k \leq 5,$$ by application of the function $f$ from \autoref{def:Zobrazeni_f}. 
For $k\in\{0,1,2,3,4\}$, the bispecial factor $F_k$ may occur in a~bispecial triplet  $(\{a,b\}, w, \{c,d\})$, where $\tt{6} \notin \{a,b\}$ or $\tt{6} \notin \{c,d\}$. Without loss of generality we assume that $a<b$ and $c<d$.  
Thanks to the previous formula from \autoref{rem:LcsLcp} for  $\lcs\{\pi(i), \pi(j)\}$, we can replace $\{a,b\}$ by $\{a,\tt{6}\}$ except for the case when $\{a,b\} = \{\tt{1},\tt{2}\}$. Analogously for the pair $\{c,d\}$ of the right extensions.  
The only bispecial factor associated with a~bispecial triplet $(\{a,b\}, F_k,\{c,d\})$ such that $\{a,b\} = \{\tt{1},\tt{2}\}$ or $\{c,d\} = \{\tt{1},\tt{2}\}$ is $F_1=\tt{0}$.  Therefore, in addition to the bispecial triplets generated from the initial triplets, we also compute the $\pi$-offsprings of the bispecial triplets $
(\{\tt{1},\tt{2}\}, \tt{0}, \{\tt{1},\tt{2}\})$ and $(\{\tt{1},\tt{6}\}, \tt{0}, \{\tt{1},\tt{2}\})$.

\begin{example}\label{ex:short} Let $w \in \{F_0,F_1,F_2\} $.   First, we find for $w$ the bispecial triplets it is  associated with and then we determine the $\pi$-offsprings.  For some of the  $\pi$-offsprings we also specify the shortest return word~$z$ to convince the reader that $|v|/|z| \leq 1$  for each $\pi$-offspring $v$  and  each return word $z$ to $v$. Thanks to \autoref{obs:BSsymmetry}, if an  $S$-image of a~$\pi$-offspring $v$ has already been examined, we do not need to care about the ratio $|v|/|z|$.  We also use the fact that if $\tilde{z}$  is the shortest return word to $\tilde{v}$ and $\tilde{v}$ is  a factor of $v$, then $|z|\geq |\tilde{z}|$ for every return word $z$ to $v$.  

\begin{itemize}
    \item $w = F_0=\varepsilon$ is associated with the bispecial triplets  $T_k$, for $k \in \{0, 1, \ldots, 5\}$.  We use \autoref{rem:LcsLcp} to get their mutually different $\pi$-offsprings 
\begin{align*}
v_0&={\tt 0}=\lcs\{\pi(\tt{0}), \pi(\tt{6})\}\lcp\{\pi(\tt{0}),\pi(\tt{6})\}, & &  z_0=v_0;\\
v_1&=\tt{02}=\lcs\{\pi(\tt{0}), \pi(\tt{6})\}\lcp\{\pi(\tt{3}),\pi(\tt{6})\}, & & z_1=v_1\tt{23223322322};\\
v_2&=\tt{022}=\lcs\{\pi(\tt{0}), \pi(\tt{6})\}\lcp\{\pi(\tt{4}),\pi(\tt{6})\}, & & |z_2|\geq |z_1|;\\
v_3&=\tt{022322}=\lcs\{\pi(\tt{0}), \pi(\tt{6})\}\lcp\{\pi(\tt{5}),\pi(\tt{6})\}, & & |z_3|\geq|z_1|.
\end{align*}

\item $w=F_1=\tt{0}$ is associated with the bispecial triplets $f(T_k)$ for $k \in \{0,1,2,3,4\}$, along with the bispecial triplets $
(\{1,2)\}, \tt{0}, \{1,2\})$ and $(\{1,6\}, \tt{0}, \{1,2\})$, see the explanation above. 
By \autoref{rem:LcsLcp}, their mutually different $\pi$-offsprings are
\begin{align*}
v_4&={\tt{00}}=\lcs\{\pi(\tt{1}), \pi(\tt{6})\}\,\pi(w)\lcp\{\pi(\tt{1}),\pi(\tt{6})\}, & & z_4=v_4\tt 1;\\
v_5&={\tt{1001}}=\lcs\{\pi(\tt{1}), \pi(\tt{2})\}\,\pi(w)\lcp\{\pi(\tt{1}),\pi(\tt{2})\}, & & z_5=v_5;\\
v_6&={\tt{001}}=\lcs\{\pi(\tt{1}), \pi(\tt{6})\}\,\pi(w)\lcp\{\pi(\tt{1}),\pi(\tt{2})\}, & & z_6=v_6;\\
v_7&={\tt{002}}=\lcs\{\pi(\tt{1}), \pi(\tt{6})\}\,\pi(w)\lcp\{\pi(\tt{3}),\pi(\tt{6})\}=S(v_2); \\
v_8&={\tt{0022}}=\lcs\{\pi(\tt{1}), \pi(\tt{6})\}\,\pi(w)\lcp\{\pi(\tt{4}),\pi(\tt{6})\}, & & |z_8|\geq |z_1|;\\
v_9&={\tt{0022322}}=\lcs\{\pi(\tt{1}), \pi(\tt{6})\}\,\pi(w)\lcp\{\pi(\tt{5}),\pi(\tt{6})\}, & & |z_{9}|\geq |z_1|.
\end{align*}

\item $w=F_2=\varphi(\tt{0})\tt{0}=\tt{010}$ is associated with the bispecial triplets $f^2(T_k)$ for $k \in \{0,1,2,3\}$. Their $\pi$-offsprings are  
\begin{align*}
v_{10}&={\tt{00100}}=\lcs\{\pi(\tt{2}), \pi(\tt{6})\}\,\pi(w)\lcp\{\pi(\tt{2}),\pi(\tt{6})\}, & & z_{10}=v_{10}\tt{2};\\
v_{11}&={\tt{001002}}=\lcs\{\pi(\tt{2}), \pi(\tt{6})\}\,\pi(w)\lcp\{\pi(\tt{3}),\pi(\tt{6})\}=S(v_3); \\
v_{12}&={\tt{0010022}}=\lcs\{\pi(\tt{2}), \pi(\tt{6})\}\,\pi(w)\lcp\{\pi(\tt{4}),\pi(\tt{6})\}=S(v_9);\\
v_{13}&={\tt{0010022322}}=\lcs\{\pi(\tt{2}), \pi(\tt{6})\}\,\pi(w)\lcp\{\pi(\tt{5}),\pi(\tt{6})\}, & & |z_{13}|\geq |z_1|. \qedhere
\end{align*}

\end{itemize}
    
\end{example}

On the other hand, every bispecial factor of $\zz$ containing $\tt 0$ is the $\pi$-offspring of some bispecial triplet of $\uu$.
\begin{lemma}\label{lem:preimage}
If $v$ is a bispecial factor in $\zz$ containing the factor $\tt 0$, then there exists a~bispecial triplet $T=(\{a,b\}, w, \{c,d\})$ in $\uu$ such that $v$ is the $\pi$-offspring of $T$.  Further, the bispecial factor $w$ associated with the bispecial triplet $T$ is uniquely determined by $v$, and if $w$ has exactly two left and two right extensions in $\uu$, then for every return word $z$ to $v$ in $\zz$, there exists a return word $r$ to $w$ in $\uu$ such that $z$ is conjugate to $\pi(r)$. 
\end{lemma}

\begin{proof} 
Let $v$ be a bispecial factor in $\zz$ containing $\tt{0}$.  Identify the first and the last occurrence of $\tt 0$ in $v$, i.e., write $v=x{\tt 0}\cdots {\tt 0}y$, where $x$ and $y$ do not contain $\tt{0}$.  (Note that the first and the last occurrence of $\tt{0}$ could coincide, in which case we have $v=x\tt{0}y$.)  By the form of $\pi$, it follows that $v=x\pi(w){\tt 0}y$, where 
\begin{itemize}
    \item $w$ is a uniquely determined (and possibly empty) bispecial factor in $\uu$, and
    \item $x=\lcs\{\pi(a), \pi(b)\}$ and ${\tt 0}y=\lcp\{\pi(c),\pi(d)\}$ for some bispecial triplet $(\{a,b\}, w, \{c,d\})$ in $\uu$.
\end{itemize}
To prove the last part of the statement, assume that $a$ and $b$ are the only left extensions of $w$ in $\uu$, and that $c$ and $d$ are the only right extensions of $w$ in $\uu$.
Let us consider a~return word $z$ to $v$ in $\zz$. Then $zv$ contains $v$ only as a prefix and a suffix, hence $zv=x\pi(w){\tt 0}y\cdots x\pi(w){\tt 0}y$. Then necessarily $zv=x\pi(rw){\tt 0}y=x\pi(r)x^{-1}x\pi(w){\tt 0}y$, where $r$ is a return word to $w$ in $\uu$. Otherwise, if $r$ is not a return word to $w$, then $zv$ contains another occurrence of $\pi(w)$, but since $w$ has only two left and two right extensions, $zv$ contains another occurrence of $x\pi(w){\tt 0}y=v$, which is a contradiction. Since $z=x\pi(r)x^{-1}$, the word $z$ is conjugate to $\pi(r)$.
\end{proof}

\begin{remark}
The only bispecial factors of $\zz$ that contain neither $\tt 0$ nor $\tt 2$ are $\tt 1$ and $\tt 3$, and their shortest return words are $\tt 1$ and $\tt 3$, respectively. 
Every other bispecial factor in $\zz$ either contains $\tt 0$ and is the $\pi$-offspring of some bispecial triplet in $\uu$, as described in \autoref{lem:preimage}, or contains $\tt 2$ and is the $S$-image of a bispecial factor containing $\tt 0$.
\end{remark}

To determine the critical exponent  of ${\bf z}$ using \autoref{thm:FormulaForE}, we must examine the bispecial factors of $\zz$ and their shortest return words. The following lemmas will limit the list of bispecial factors that need to be examined.

\begin{lemma}\label{coro:nonpalindromic}
Let $v$ be a bispecial factor in $\zz$ that is the $\pi$-offspring of a~bispecial triplet associated with a nonpalindromic bispecial factor $w$ of $\uu$. 
Then $\frac{|v|}{|z|}<1$ for every return word $z$ to~$v$.
\end{lemma}
\begin{proof}
Every nonpalindromic bispecial factor $w$ in $\uu$ has exactly two left and two right extensions. 
Consequently, by \autoref{lem:preimage}, a~return word $z$ to $v$ is conjugate to $\pi(r)$ for some return word $r$ to $w$, hence $|z| = |\pi(r)|$. By \cite[Lemma 10]{DvPe_LimRT2025}, the return word $r$ to $w$ satisfies $\vv{r}\geq \vv{w}+\vv*{e}{6}\ \text{or} \ \vv{r}\geq \vv{w}+\vv*{e}{5} + \vv*{e}{4}+\vv*{e}{3}+\vv*{e}{2}+\vv*{e}{1}+\vv*{e}{0}.$ Hence we have
\[
|z|\geq |\pi(w)|+|\pi(\tt{6})|=|\pi(w)|+13
\ \ \text{or} \ \ 
|z|\geq |\pi(w)|+\text{$\sum_{k=0}^5$}|\pi({k})|=|\pi(w)|+17\,.
\]
But since $v$ is the $\pi$-offspring of some bispecial triplet $(\{a,b\},w,\{c,d\})$, using \autoref{obs:poProjekci} and the form of $\pi$, we have $|v|\leq |\pi(w)|+2|\pi(\tt{5})|\leq |\pi(w)|+12$.  Altogether, we have
\[
\frac{|v|}{|z|}\leq \frac{|\pi(w)|+12}{|\pi(w)|+13}<1\,.\qedhere
\]
\end{proof}

\begin{example}\label{ex:Simages} 
Let us show that $\pi$-offsprings arising from   palindromic and nonpalindromic bispecial factors in $\uu$ can be images of one another under the antimorphism $S$.  This is the case for the $\pi$-offspring of the following pairs.  In the calculations below, we use the equality $F_kk = \varphi^k(\tt{0})$ and \autoref{ex:triplets}.

\begin{enumerate}

\item $f^3(T_1)$ and  $f(T_5)$
    
The $\pi$-offspring  of  $f^3(T_1) =(\{\tt{3},\tt{6}\}, F_3, \{\tt{4},\tt{6}\})$  is  
\[
v=\tt{2}  \pi(F_3)\tt{022}= \tt{2}    \pi(F_3\tt{3})\tt{2} = \tt{2}  \pi(\varphi^3(\tt{0}))\tt{2},
\]
and the $\pi$-offspring  of  $f(T_5) =(\{\tt{1},\tt{6}\}, \tt{06}, \{\tt{0},\tt{6}\})$ is $u=\pi(\tt{06})\tt{0} = \tt{0}\pi(\tt{6})\tt{0}$.  Using \eqref{eq:obrazByS} we see that $v=S(u)$. 

\item  $f^4(T_1)$ and   $f^2(T_5)$

The $\pi$-offspring  of       $f^4(T_1) = (\{\tt{4},\tt{6}\}, F_4,  \{\tt{5},\tt{6}\})$   is 
\[
v=\tt{22}  \pi(F_4)\tt{022322}= \tt{22}   \pi(F_4\tt{4})\tt{322} = \tt{22}  \pi(\varphi^4(\tt{0}))\tt{322},
\]
and the $\pi$-offspring  of  $f^2(T_5) =(\{\tt{2},\tt{6}\}, \varphi(\tt{06})\tt{0}, \{\tt{1},\tt{6}\})$ is 
\[
u=\pi(\varphi(\tt{06})\tt{0})\tt{0} = \tt{001}\pi(\varphi(\tt{6}))\tt{00},
\]
hence $v=S(u)$. 

\item  $f^3(T_2)$ and $f^2(T_4)$

The $\pi$-offspring  of       $f^3(T_2) = (\{\tt{3},\tt{6}\}, F_3,  \{\tt{5},\tt{6}\})$   is
\[
\tt{2}  \pi(F_3)\tt{022322}= \tt{2}    \pi(F_3\tt{3})\tt{2322} = \tt{2}  \pi(\varphi^3(\tt{0}))\tt{2322},
\]
and the $\pi$-offspring  of  $f^2(T_4) =(\{\tt{2},\tt{6}\}, \varphi(\tt{0})\tt{06}, \{\tt{0},\tt{6}\})$ is 
\[
u=\pi(\varphi(\tt{0})\tt{0 6}))\tt{0} = \tt{0010}\pi(\tt{6})\tt{0},
\]
hence $v=S(u)$.  \qedhere
\end{enumerate}
\end{example}

\begin{lemma}\label{lem:T6only}  The critical exponent $E(\zz) $ of $\zz$ equals 
\[
 1 + \sup\left\{\tfrac{|v|}{|z|} \,\colon\,  v \text{ is the $\pi$-offspring of $f^n(T_6)$ for some $n\in\N$ and $z$ is a return word to $ v$ in $\zz$ }\right\}.
\]
\end{lemma}
\begin{proof} By \autoref{thm:FormulaForE}, the critical exponent of $\zz$ is the maximum ratio $|v|/|z|$, where $v$ is a bispecial  factor of $\zz$ and $z$ is a return word to $v$ in $\zz$.  As $E(\xx) > 2$ for every  rich infinite word $\xx$, it is enough to consider the cases where $|v|/|z|>1$.  The closedness of the language of $\zz$ under the antimorphism  $S$ implies that $v$ is a bispecial factor of $\zz$ and $zv$ a~complete return word to $v$ in $\zz$ if and only if the same is true for $S(v)$ and $S(zv)$, respectively. Obviously, $|v|/|z| = |S(v)|/|S(z)|$.    Therefore, to calculate the critical exponent, it is sufficient to consider one of the pair of bispecial factors $v$  and $S(v)$. 

By  \autoref{coro:nonpalindromic}, we can  focus on  bispecial factors   which are the   $\pi$-offsprings of bispecial triplets associated with a~palindromic bispecial factor of $\uu$. 

Palindromic bispecial factors of $\uu$ are associated with the bispecial triplets $f^n(T_6)$, $f^n(T_0)$, for $n \in \N$, and with the bispecial triplets  $f^n(T_k)$ for   $k \in \{1,2,3,4,5\}$  and  $n \leq 5-k$.   

By \autoref{ex:short}, we do not need to discuss the cases $f^n(T_k)$ for $k =0,1,\ldots, 5$ and $n \leq 2$,  as the associated bispecial factors   belong  to  $\{F_0, F_1,F_2\}$. 

By   \autoref{ex:Simages},  the $\pi$-offsprings of $f^3(T_1)$, $f^4(T_1)$, and $f^3(T_2)$ are images under $S$ of $\pi$-offsprings of bispecial triplets associated with a nonpalindromic bispecial factor of $\uu$. Therefore, they can be excluded from our consideration, too. 

It remains to discuss the $\pi$-offsprings of bispecial triplets of the form $f^n(T_0)$ with $n\geq 3$. 
For every $n \in \N$, we denote $v_n$ the  $\pi$-offspring of the bispecial triplet $f^{n+3}(T_0)$.   Note that $S(v_n)$ is a bispecial factor in $\zz$ and thus it is the  $\pi$-offspring of some bispecial triplet. As $v_n$ is a palindrome, $S(v_n)\neq v_n$.  Since  the length  $|v_n|$ is increasing with $n$, 
$S(v_n)$ cannot coincide with $v_m$ for some $m \in \N$.  Hence,  $S(v_n)$ is either  the $\pi$-offspring of a~bispecial triplet $f^m(T_6)$ for some $m \in \N$ or  $S(v_n)$ is the $\pi$-offspring  of a~bispecial triplet, which according to the previous discussion does not need to be considered.

We can therefore conclude that to determine the critical exponent $E(\zz)$ it is sufficient to consider $\pi$-offsprings of  $f^n(T_6)$ as stated in the lemma.  
\end{proof}

\subsection{Proof of \autoref{thm:main}}
Thanks to \autoref{lem:T6only}, we need to study only 
bispecial factors arising by projection from the triplets $f^n(T_6)$ and their shortest return words. 
For all $n\in\N$, we write $f^n(T_6)= (\{a_n,\tt{6}\}, w_n, \{a_n,\tt{6}\})$, and we let $r_n$ denote the shortest return word to the bispecial factor $w_n$ in $\zz$. In particular, $a_0=\tt{0}$ and $w_0=r_0=\tt{6}$.  The left extension and the right extension of the bispecial factor $w_n$ consists of the letters $a_n$ and $\tt 6$, where $a_{n} = (a_{n-1}+1)\bmod 6 $. As $a_0 = \tt 0$, we can write $a_n =n\bmod 6$.  

The incidence matrices of the morphisms $\varphi$ and $\pi$ are given by
\begin{equation}\label{eq:matrices}
M=M_{\varphi}=\left(\begin{smallmatrix}
1&1&1&1&1&1&1\\
1&0&0&0&0&0&0\\
0&1&0&0&0&0&0\\
0&0&1&0&0&0&0\\
0&0&0&1&0&0&0\\
0&0&0&0&1&0&0\\
0&0&0&0&0&1&2
\end{smallmatrix}\right) \ \ \ \text{ and } \ \ \   M_{\pi}=\left(\begin{smallmatrix}
1&1&1&1&1&1&1\\
0&1&2&0&0&0&0\\
0&0&0&1&2&4&8\\
0&0&0&0&0&1&4
\end{smallmatrix}\right)\,.
\end{equation}
By \cite[Proposition 15]{DvorakovaKloudaPelantova2025}, we have
\begin{equation}\label{eq:BSze6}  
\vv*{w}{n} = 
\begin{cases}
    M\vv*{w}{n-1} + \vv*{e}{0} + 2\vv*{e}{6}, & \text{if } n\bmod 6=0\,;\\
M\vv*{w}{n-1} + \vv*{e}{0},  & \text{otherwise\,; }
\end{cases}
\qquad \text{and }\qquad \vv*{r}{n} =  M\, \vv*{r}{n-1}\,.
    \end{equation}
Note that the Parikh vectors $\vv*{w}{n}$  and $\vv*{r}{n}$ belong to $\mathbb{R}^7$. 

\begin{example}\label{ex:firstRatioInU} Let us apply  \eqref{eq:BSze6}  to compute $\vv*{w}{i}$ for $i = 0,1,\ldots, 5$.  As $\vv*{w}{0} = \vv*{e}{6}$, we obtain 
$$
\vv*{w}{1} = \left(\begin{smallmatrix}    2   \\
     0   \\
     0   \\
     0   \\
     0   \\
     0\\
     2\end{smallmatrix}\right)
 \ \vv*{w}{2} = \left(\begin{smallmatrix}  
    5   \\
     2   \\
     0   \\
     0   \\
     0   \\
     0\\
     4
\end{smallmatrix}\right),\ \vv*{w}{3} =\left(\begin{smallmatrix}  
    12  \\
     5   \\
     2   \\
     0   \\
     0   \\
     0\\
     8
\end{smallmatrix}\right),\ \vv*{w}{4} = \left(\begin{smallmatrix}  
    28   \\
     12   \\
     5  \\
     2   \\
     0   \\
     0\\
     16
\end{smallmatrix}\right), \ \vv*{w}{5} = \left(\begin{smallmatrix}  
    64   \\
     28   \\
     12  \\
     5  \\
     2   \\
     0\\
     32
\end{smallmatrix}\right).
$$
The lengths of $w_i$ thus satisfy
$$
|w_0|=1, \ |w_1|=4, \ |w_2|=11, \ |w_3|=27, \ |w_4|=63, \  |w_5|=143. 
$$
Analogously, we get

$$
\vv*{r}{1} = \left(\begin{smallmatrix}  
    1  \\
     0   \\
     0   \\
     0   \\
     0   \\
     0\\
     2
\end{smallmatrix}\right),  \ \vv*{r}{2} = \left(\begin{smallmatrix}  
    3   \\
     1   \\
     0   \\
     0   \\
     0   \\
     0\\
     4
\end{smallmatrix}\right),\ \vv*{r}{3} = \left(\begin{smallmatrix}  
    8  \\
     3   \\
     1   \\
     0   \\
     0   \\
     0\\8
\end{smallmatrix}\right),\ \vv*{r}{4} = \left(\begin{smallmatrix}  
    20   \\
     8   \\
     3  \\
     1   \\
     0   \\
     0\\16
\end{smallmatrix}\right), \ \vv*{r}{5} = \left(\begin{smallmatrix}  
    48  \\
     20   \\
     8  \\
     3  \\
     1   \\
     0\\32
\end{smallmatrix}\right),
$$
and 
$$
|r_0|=1, \ |r_1|=3, \ |r_2|=8, \ |r_3|=20, \ |r_4|=48,\ |r_5|=112. 
$$
Note that by \autoref{thm:FormulaForE} 
$$E(\uu_7)\geq 1+ \frac{|w_2|}{|r_2|} = 1+\frac{11}{8} = 2. 375 > E^*(\uu_7) \doteq 2.1197196907\,. $$
It can be shown that $E(\uu_7)= 2.375$. 
\end{example}

The following lemma follows directly from \eqref{eq:BSze6} and the simple fact that $\vv*{e}{0} = (M-2I)\vv*{e}{6}$.

\begin{lemma}\label{lem:BSExplicit} Let $ N \in \N$ and $i \in \{0,1,\ldots, 5\} $, and let $n = 6N+i$. Then $\vv*{r}{n} = M^n\,\vv*{e}{6}$  and 
$$
\vv*{w}{n} = \vv*{r}{n} + \sum_{j=0}^{n-1} M^j\,\vv*{e}{0} +  2M^i\sum_{j=0}^{N-1} M^{6j}\,\vv*{e}{6} =  \vv*{r}{n}  +  \sum_{j=0}^{n-1} M^j(M-2I)\,\vv*{e}{6}  +  2M^i\sum_{j=0}^{N-1} M^{6j}\,\vv*{e}{6}\,.
$$
\end{lemma}

\bigskip
Let us now focus on bispecial factors and their shortest return words in the word $\zz = \pi(\uu)$. Let $v_n$ denote the bispecial factor of $\zz$ arising from $w_n$ as described in \autoref{lem:preimage}.  In other words, $v_n$ is the $\pi$-offspring of $f^n(T_6)$. Let $z_n$ denote the shortest return word to $v_n$ in $\zz$. As $a_n = n \bmod 6$, the form of the projection $\pi$ and \autoref{rem:LcsLcp} imply 
\[
\tt{0}\lcs\{\pi(a_n), \pi(\tt{6})\} = \lcp\{\pi(a_n), \pi(\tt{6})\} = 
\begin{cases}\pi(a_n), & \text{if } n\bmod 6\in\{3,4,5\}\,;\\
\pi(\tt{0}), & \text{otherwise.}
\end{cases}
\]

To determine the lengths of $v_n$ and $z_n$ needed for the evaluation of the critical exponent of $\zz$, we apply \autoref{lem:preimage} to deduce that for every $n \in \N$, we have
\begin{equation}\label{eq:delky}
|v_n| = |\pi(w_n)| + {\rm  add}_n \qquad \text{and} \qquad |z_n| = |\pi(r_n)|, \quad \text{where} \end{equation}
\begin{equation}\label{eq:plusSomthing}
{\rm  add}_n  = 
\begin{cases}
    1,  & \text{if } n\bmod 6 \in\{0,1,2\}\,;\\
    3,  & \text{if } n=3\bmod 6\,;\\
    5,  & \text{if } n=4\bmod 6\,;\\
    11, & \text{if } n=5\bmod 6\,.\\
\end{cases}
    \end{equation}

To compute the length of factors of $\zz$, we will use the Parikh vectors and the matrices of the corresponding morphisms.  
\begin{observation} If $v \in \Sigma_4^*$, then its length is $ |v|= (1,1,1,1)\vv{v}$. If moreover $v = \pi (w)$, where $w \in \Sigma_7^*$, then
 \begin{equation}\label{eq:DELKA}
     |\pi(w)| = (1,1,1,1)\vv{\pi(w)} = (1,1,1,1)\,M_\pi\, \vv{w} = (1,2,3,2,3,6,13)\, \vv{w}\,.
 \end{equation}   
\end{observation}
We denote 
\begin{equation}\label{h}
\vv{h}=(1,1,1,1)\,M_\pi=(1,2,3,2,3,6,13)\,.
\end{equation}
\begin{example}\label{ex:firstRatio} Using \autoref{ex:firstRatioInU} along with \eqref{eq:delky} and \eqref{eq:DELKA}, we obtain the information in the table below.
\[
\begin{array}{c|c|c|c|c}
i & \text{length of }\pi(w_i)& \text{length of }v_i&  \text{length of } z_i& 1+\tfrac{|v_i|}{|z_i|} \rule[-2ex]{0pt}{2ex}\\
\hline
0&13&14&13& 2.07692\\
1&28&29&27&2.07407 \\
2&61&62&57&2.08772 \\
3&132&135&121&2.11570 \\
4&279&284&255&2.11373 \\
5&588&599&537& 2.11546
\end{array}  \qedhere
\]
\end{example}
The previous example illustrates that the values $1+\tfrac{|v_i|}{|z_i|}
$ are smaller than $E^*(\uu_7)$ for $i<6$. 
Our goal is to show that this is the case for all $i \in \N$. 
In order to express the lengths of $v_n$ and $z_n$ in a simple way, we will use the results from the appendix of~\cite{DvPe_LimRT2025}.

 \begin{lemma}\label{lem:reseniRekurence6}  Let $\Lambda=\mu_1,  \mu_2,\mu_3, \mu_4$ denote the roots of the polynomial $t^4-2t^3-1$. Then 
 $$\vv{h} M^n \vv*{e}{6}=\sum_{j=1}^4 c_j \mu_j^n\,,\quad \text{where}\ c_j= \frac{(\mu_j^3+\mu_j+2)\mu_j(\mu_j-1)}{4\mu_j-6}\,.$$
\end{lemma}
\begin{proof}
We can draw the following information from the appendix of~\cite{DvPe_LimRT2025}:
\begin{itemize}
\item The characteristic polynomial of the matrix $M$ is given by
$$\chi(t)=(t^4-2t^3-1)(t^3-t^2-t-1)=\tfrac{1}{t-1}\Bigl(t^{6}\bigl(t-2)^2-1\Bigr)\,,$$
where all eigenvalues in the spectrum of $M$ are distinct and the spectrum contains just two eigenvalues outside the closure of the unit circle, namely 
$\Lambda \in (2,3)$ and $\beta \in (1,2)$. They  satisfy  $\Lambda^{4}-2\Lambda^3=1$ and $\beta^{4}-2\beta^3=-1.$
Let $S_\Lambda$ denote the set of algebraic conjugates of $\Lambda$, and let $S_\beta$ denote the set of algebraic conjugates of $\beta$.
Thus the spectrum of $M$ is $S=S_\Lambda \cup S_\beta$.
\item The right eigenvector of $M$ corresponding to an eigenvalue $\lambda \in S$ is 
\begin{equation}\label{eq:eigenvector}
\vv*{v}{\lambda} = \bigl(\lambda^{5}(\lambda-2),\lambda^{4}(\lambda-2), \lambda^{3}(\lambda-2), \lambda^{2}(\lambda-2), \lambda(\lambda-2), \lambda-2, 1\bigr)^T;
\end{equation}
the left row eigenvector to $\lambda$ is 
$$\vv*{u}{\lambda} = \bigl(1, \lambda -1, \lambda^2 - \lambda -1, \ldots, \lambda^{6} - \lambda^{5} - \lambda^4-\lambda^3-\lambda^2 - \lambda -1\bigr)\,.$$
In particular,  if $\lambda \in S_{\Lambda}$, then the last  entry of  $\vv*{u}{\lambda}$ is 
$$\frac{1}{\lambda -1}\bigl(\lambda^{7} - 2 \lambda^{6} +1\bigr) = \lambda^3\,.$$
\end{itemize}
Setting $\vv{g}=\vv*{e}{6}$ in Lemma 20 from the appendix of~\cite{DvPe_LimRT2025}, we get 
$$
\vv{h}M^n\vv*{e}{6}=\sum_{\lambda \in S} c_\lambda \lambda^n\,, \quad \text{where}\ c_\lambda= \frac{1}{\vv*{u}{\lambda}\,\vv*{v}{\lambda}}\, \bigl(\vv{h}\,\vv*{v}{\lambda}\bigr) \, \bigl(\vv*{u}{\lambda}\, \vv*{e}{6}\bigr)\,.
$$
It is easy to check that $\vv{h}\,\vv*{v}{\lambda}=0$ for $\lambda \in S_\beta$ and $c_\lambda=\frac{(\lambda^3+\lambda+2)\lambda(\lambda-1)}{4\lambda-6}$ for $\lambda \in S_\Lambda$. The statement follows.
\end{proof}




Let us list approximate values of the roots $\mu_i$ and the coefficients $c_i$ that will be used later.
$$\begin{array}{rcl}
\mu_1&\doteq &2.1069134037622186,\\
\mu_2&\doteq &-0.71667249,\\
\mu_3&\doteq & 0.3048767+i 0.754529173\\
\mu_4& = &\overline{\mu_3}
\end{array} \qquad \text{and}\qquad \begin{array}{rcl}
c_1&\doteq & 12.9304\\
c_2&\doteq & -0.127\\
c_3&\doteq & 0.0983+i 0.2612\\
c_4&=&\overline{c_3}
\end{array}$$
\begin{proposition}\label{prop:BSformT6} For $n\in \N$, let $v_n$ denote the $\pi$-offspring of $f^n(T_6)$, and let $z_n$ denote the shortest return word to $v_n$ in $\zz$.  Then for all $n\in\N$, we have
$$1+\frac{|v_n|}{|z_n|} \leq E^*(\uu_7)\,.$$
\end{proposition}
\begin{proof} 
Let $n =6N+i$, where $N\in\N$ and $i\in\{0,1,\ldots,5\}$. We consider only the case that $N\geq 1$, as for $n <6$ the validity of the statement is demonstrated in \autoref{ex:firstRatio}. To compute the ratio ${|v_n|}/{|z_n|}$, we use~\eqref{eq:delky} and~\eqref{eq:plusSomthing}. 

Since $E^*(\uu_7)= 1+\frac{1}{3-\Lambda}$, we need to show that $\frac{|v_n| - |z_n|}{|z_n|} \leq \frac{\Lambda-2}{3-\Lambda}$. Using \autoref{lem:BSExplicit} along with \eqref{eq:delky}, we can rewrite this inequality as
\begin{equation}\label{eq:wantedIneq}
\vv{h}\,A\,\vv*{e}{6} + {\rm add}_n\  \leq \ \tfrac{\Lambda-2}{3-\Lambda} \,\vv{h}\, M^n\, \vv*{e}{6}, 
\end{equation}
where
$$
A = (M-2I)\sum_{j=0}^{n-1}M^j +2 M^i \sum_{j=0}^{N-1}M^{6j}\,. 
$$
Since every matrix $H$ without eigenvalue $1$ satisfies $\sum_{j=1}^{k-1}H^j = (H-I)^{-1}(H^k-I)$, we can write 
$$
A = \underbrace{\Bigl((M-I)^{-1}(M-2I)+2(M^6-I)^{-1}\Bigr)}_{A_1}M^n - \underbrace{\Bigl((M-I)^{-1}(M-2I) +2(M^6-I)^{-1}M^i\Bigr)}_{A_2}\,. 
$$

To confirm \eqref{eq:wantedIneq}, we will prove the following two inequalities.
$$
\begin{array}{ccccccc}
 \text{Inequality 1:}\qquad &  \vv{h}\,A_1\,M^n \vv*{e}{6} &&&\leq &\, \tfrac{\Lambda-2}{3-\Lambda}\,\vv{h}\,M^n\vv*{e}{6}   &+ \frac{6}{13}
   \\
 \text{Inequality 2:} \qquad & & - \vv{h}\,A_2\,\vv*{e}{6} &
+\,{\rm add}_n &\leq & &- \frac{6}{13}
\end{array}
$$
Summing these two inequalities gives \eqref{eq:wantedIneq}. 

\begin{description}
\item[Inequality 1.]
The Cayley-Hamilton theorem says that $M^6 (M-2I)^2 = I$. This implies that 
\[
(M^6 - I)(M-2I)^2 = I - (M-2I)^2 = (3I-M)(M-I)\,,
\]
and it follows that 
\[
(M^6 - I)^{-1}=(M-I)^{-1}(3I-M)^{-1}(M-2I)^2\,.
\]
Consequently,
\begin{align*}
    A_1 &=  (M-I)^{-1}(M-2I)+2(M-2I)^2(3I-M)^{-1}(M-I)^{-1} \\
    &= (M-2I)(3I-M)^{-1}\,.
\end{align*}
We apply \autoref{lem:reseniRekurence6} to rewrite  Inequality 1. So we need to demonstrate that  
$$\sum_{j=2}^4\mu_j^{n} c_j\left(\tfrac{\mu_j -2}{3-\mu_j}-\tfrac{\Lambda-2}{3-\Lambda}\right)\leq \tfrac{6}{13}\,.$$
We exploit the explicit values of $\mu_j$ and $c_j$  and the fact that $|\mu_j| < 1$ for $j \in \{2,3,4\}$. For $n\geq 6$, we find that 
$$\Bigl|\sum_{j=2}^4\mu_j^{n} c_j\Bigl(\tfrac{\mu_j -2}{3-\mu_j}-\tfrac{\Lambda-2}{3-\Lambda}\Bigr)\Bigr|\leq \sum_{j=2}^4 |\mu_j|^6 \cdot |c_j|\,\Bigl(\Bigl|\tfrac{\mu_j -2}{3-\mu_j}-\tfrac{\Lambda-2}{3-\Lambda}\Bigr|\Bigr) \ \ \doteq \ 0,142 \  \leq \  \tfrac{6}{13}\,,$$
as desired.

\item[Inequality 2.]  
To verify Inequality 2, we recall the values of ${\rm add}_i$ listed in \eqref{eq:plusSomthing}, and  calculate the values of the expression $\vv{h}\,A_2\,\vv*{e}{6}- \tfrac6{13}$ in the table below.
$$\begin{array}{c|c|c|c|c|c|c}
     i&0&1&2&3&4&5 \\
     \hline

     {\rm add}_i &1&1&1&3&5&11 \rule[-1ex]{0pt}{1ex}\\      
     \vv{h}\,A_2\,\vv*{e}{6}  - \tfrac6{13}& \phantom{11} 1\phantom{11}  & 1+ \tfrac{12}{13} &2+ \tfrac{9}{13}& \phantom{11} 3 \phantom{11} & 5+ \tfrac{12}{13}&12+ \tfrac{9}{13}\rule[-1ex]{0pt}{1ex} 
\end{array}
$$
Consequently, Inequality 2 is valid as well.  \qedhere 
\end{description}    
\end{proof}

\begin{proof}[Proof of \autoref{thm:main}]
On the one hand, in \autoref{sec:BS}, we have shown that for each bispecial factor $v$ of $\zz=\pi(\uu_7)$ and any of its return words $z$ the inequality $1+\frac{|v|}{|z|}\leq E^*(\uu_7)$ holds. Thus by \autoref{thm:FormulaForE}, we have $E(\zz)\leq E^*(\uu_7)$.
On the other hand, by~\cite[Theorem 9]{Ochem}, we have the relation $E^*(\pi(\uu_7))=E^*(\uu_7)$. Since $E(\zz)\geq E^*(\zz)=E^*(\uu_7)$, we obtain the equality $E(\zz)=E^*(\uu_7)$, as desired. 
\end{proof}

\section{Comments and open questions}
\label{Section:OpenQuestions}

\begin{enumerate}

    \item A first step towards proving \autoref{con:opet} for every $d\geq 5$ would be to find a projection $\pi=\pi_d:\Sigma^*_{2d-1}\to \Sigma^*_d $ such that $\pi({\bf u}_{2d-1})$ is rich and    $E(\pi({\bf u}_{2d-1}))=E^*({\bf u}_{2d-1})$.  This would prove the inequality $\RT(\mathcal{C}_{d}) \leq  E^*(\uu_{2d-1})$. 

    \item  We explained that to determine the repetition threshold $RT(\mathcal{C}_d)$, it is enough to consider  uniformly recurrent words, since the threshold is equal to the infimum of the set $$S_d=\{E(\uu): \uu \text{\ is a uniformly recurrent rich word over $\Sigma_d$}\}\,.$$ 
    In fact,  for $d=2,3,4$,  the value $\RT(\mathcal{C}_d)$ is the minimum of $S_d$, and moreover, it is an isolated point of $S_d$. What is the second smallest element of $S_d$? 
    Or at least, for $d\geq 5$, how can we find a rational number $t_d$ with a small denominator such that $t_d$ is between $\RT(\mathcal{C}_d)$ and $\inf(S_d\backslash\{\RT(\mathcal{C}_d)\})$?  In other words, what should the rational number $t_d$ be in \autoref{con:structure}?
    
    \item Roughly speaking, we have observed that applying the morphism $\varphi=\varphi_d$ repeatedly to a bookended repetition tends to make the repetition ``more top-heavy''.  As evidence, we note that every $t$-persistent bookended repetition $R$ that we encountered in the proofs of \autoref{Theorem:TernaryStructure} and \autoref{Theorem:QuaternaryStructure} had a tipping point of at most $69$, i.e., we obtained a top-heavy $t$-power after applying $\varphi$ to $R$ at most $69$ times.   
    We suspect that for the values of $t$ and the morphisms $\varphi$ and $\hat{\pi}$ (or $\hat{\pi}_i$) used in this paper, every $t$-persistent bookended repetition has a tipping point, but we do not have a proof.  
    For $d\geq 2$, are there simple conditions on $t$ and $\rho:\Sigma_d^*\rightarrow \Sigma_{d'}^*$ that guarantee that every $t$-persistent bookended repetition with respect to $\varphi_d$ and $\rho$ has a tipping point?
    
\end{enumerate}

\subsection*{Acknowledgements}

Research of Lucas Mol is supported by NSERC grant RGPIN-2021-04084.

\newpage

\appendix

\section{More on the backtracking searches used in the proof of \autoref{Lemma:TernaryInductiveStep}}
\label{Appendix:BacktrackingDetails}

The proof of \autoref{Lemma:TernaryInductiveStep} involves many backtracking searches.  In this appendix, we provide more detailed descriptions/illustrations of two of these searches for the curious reader.  We stress that all of the data reported here was output by our computer programs.  We have merely reformatted it here for ease of reading.  A reminder that our code (and its output) is available at the link below.
\begin{center}
\url{https://github.com/lgmol/Structure-Theorems-for-Rich-Words/}
\end{center}

\begin{example}\label{Example:ReturnsTo0}
Consider the search for all possible complete return words $w$ to $\tt{0}$ over $\Sigma_5$ such that 
\begin{enumerate}[label=(\roman*)]
    \item $\varphi^2(w)$ is rich (hence $w$ is rich), and
    \item $w$ contains no $\tfrac{16}{7}$-persistent bookended repetition with tipping point at most $100$.
\end{enumerate}   
The entire search is too large to be illustrated here, so we illustrate the search for all such return words with prefix $\tt{03}$ in \autoref{Figure:ReturnsWithPrefix03}.  Each green leaf in the tree corresponds to a possible complete return word to $\tt{0}$ in $\xx$ with prefix $\tt{03}$.  Each yellow leaf in the tree corresponds to a word that is not rich, so it fails condition (i).  Finally, each red leaf in the tree corresponds to a word that contains a $\tfrac{16}{7}$-persistent bookended repetition with tipping point at most $100$ as a suffix, so it fails condition (ii).  We provide more details concerning the red leaves below.%
\begin{figure}
        \centering
        \begin{forest}
        [
        $\tt{03}$, for tree={grow=0,l=2.5cm}
        [
            $\tt{4}$
            [
                $\tt{4}$,draw,fill=red
            ]
            [
                \phantom{$\tt{0}$},no edge
            ]
            [
                $\tt{3}$
                [
                    $\tt{4}$,draw,fill=red
                ]
                [
                    $\tt{3}$,draw,fill=red
                ]
                [
                    $\tt{2}$,draw,fill=red
                ]
                [
                    $\tt{1}$,draw,fill=red
                ]
                [
                    $\tt{0}$,draw,fill=green
                ]
            ]
            [
                $\tt{2}$
                [
                    $\tt{4}$,draw,fill=red
                ]
                [
                    $\tt{3}$,draw,fill=yellow
                ]
                [
                    $\tt{2}$,draw,fill=red
                ]
                [
                    $\tt{1}$,draw,fill=red
                ]
                [
                    $\tt{0}$,draw,fill=yellow
                ]
            ]
            [
                $\tt{1}$
                [
                    $\tt{4}$,draw,fill=red
                ]
                [
                    $\tt{3}$,draw,fill=yellow
                ]
                [
                    $\tt{2}$,draw,fill=red
                ]
                [
                    $\tt{1}$,draw,fill=red
                ]
                [
                    $\tt{0}$,draw,fill=yellow
                ]
            ]
            [
                $\tt{0}$,draw,fill=yellow
            ]
        ]
        [
            \phantom{$\tt{0}$},no edge
        ]
        [
            \phantom{$\tt{0}$},no edge
        ]
        [
            \phantom{$\tt{0}$},no edge
        ]
        [
            $\tt{3}$,draw,fill=red
        ]
        [
            $\tt{2}$
            [
                $\tt{4}$,draw,fill=red
            ]
            [
                $\tt{3}$,draw,fill=red
            ]
            [
                $\tt{2}$,draw,fill=red
            ]
            [
                $\tt{1}$,draw,fill=red
            ]
            [
                $\tt{0}$,draw,fill=yellow
            ]
        ]
        [
            $\tt{1}$
            [
                $\tt{4}$,draw,fill=red
            ]
            [
                $\tt{3}$,draw,fill=red
            ]
            [
                $\tt{2}$,draw,fill=red
            ]
            [
                $\tt{1}$,draw,fill=red
            ]
            [
                $\tt{0}$,draw,fill=yellow
            ]
        ]
        [
            \phantom{$\tt{0}$},no edge
        ]
        [
            $\tt{0}$,draw,fill=green
        ]
        ]
        \end{forest}
        \caption{An illustration of the backtracking search for complete return words $w$ to $\tt{0}$ with prefix $\tt{03}$ satisfying conditions (i) and (ii) of \autoref{Example:ReturnsTo0}.}
        \label{Figure:ReturnsWithPrefix03}
    \end{figure}
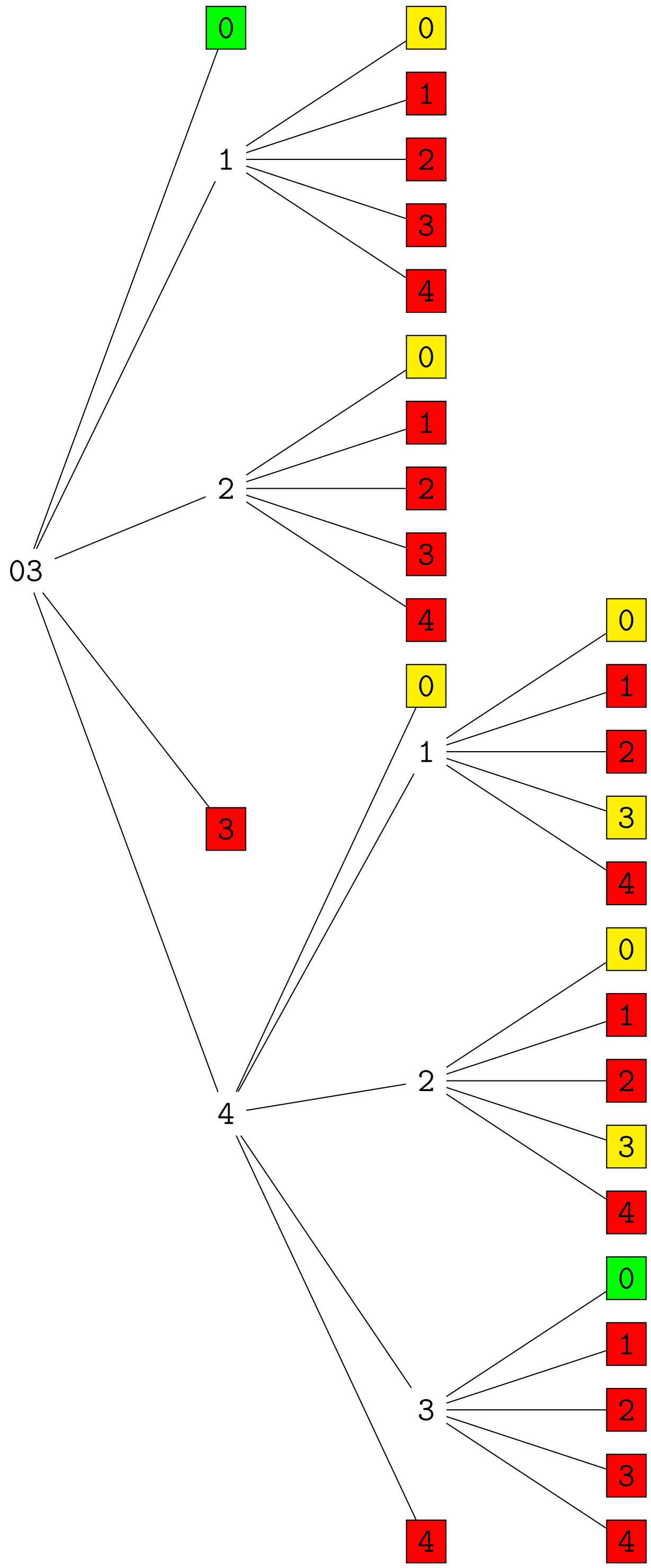%
    \begin{itemize}
        \item $\tt{0311}$ has suffix $\tt{11}$, and $(\varepsilon,\tt{11},\varepsilon,\tt{1})$ is $\tfrac{16}{7}$-persistent with tipping point $6$,
        \item $\tt{0312}$ has suffix $\tt{312}$, and $(\tt{3},\tt{1},\tt{2},\tt{1})$ is $\tfrac{16}{7}$-persistent with tipping point $5$,
        \item $\tt{0313}$ has suffix $\tt{313}$, and $(\tt{3},\tt{1},\tt{3},\tt{1})$ is $\tfrac{16}{7}$-persistent with tipping point $5$,
        \item $\tt{0314}$ has suffix $\tt{314}$, and $(\tt{3},\tt{1},\tt{4},\tt{1})$ is $\tfrac{16}{7}$-persistent with tipping point $5$,
        \item $\tt{0321}$ has suffix $\tt{321}$, and $(\tt{3},\tt{2},\tt{1},\tt{2})$ is $\tfrac{16}{7}$-persistent with tipping point $5$,
        \item $\tt{0322}$ has suffix $\tt{22}$, and $(\varepsilon,\tt{22},\varepsilon,\tt{2})$ is $\tfrac{16}{7}$-persistent with tipping point $8$,
        \item $\tt{0323}$ has suffix $\tt{323}$, and $(\tt{3},\tt{2},\tt{3},\tt{2})$ is $\tfrac{16}{7}$-persistent with tipping point $3$,
        \item $\tt{0324}$ has suffix $\tt{324}$, and $(\tt{3},\tt{2},\tt{4},\tt{2})$ is $\tfrac{16}{7}$-persistent with tipping point $3$,
        \item $\tt{033}$ has suffix $\tt{33}$, and $(\varepsilon,\tt{33},\varepsilon,\tt{3})$ is $\tfrac{16}{7}$-persistent with tipping point $14$,
        \item $\tt{03411}$ has suffix $\tt{11}$, and $(\varepsilon,\tt{11},\varepsilon,\tt{1})$ is $\tfrac{16}{7}$-persistent with tipping point $6$,
        \item $\tt{03412}$ has suffix $\tt{412}$, and $(\tt{4},\tt{1},\tt{2},\tt{1})$ is $\tfrac{16}{7}$-persistent with tipping point $5$,
        \item $\tt{03414}$ has suffix $\tt{414}$, and $(\tt{4},\tt{1},\tt{4},\tt{1})$ is $\tfrac{16}{7}$-persistent with tipping point $5$,
        \item $\tt{03421}$ has suffix $\tt{421}$, and $(\tt{4},\tt{2},\tt{1},\tt{2})$ is $\tfrac{16}{7}$-persistent with tipping point $5$,
        \item $\tt{03422}$ has suffix $\tt{22}$, and $(\varepsilon,\tt{22},\varepsilon,\tt{2})$ is $\tfrac{16}{7}$-persistent with tipping point $8$,
        \item $\tt{03424}$ has suffix $\tt{424}$, and $(\tt{4},\tt{2},\tt{4},\tt{2})$ is $\tfrac{16}{7}$-persistent with tipping point $3$,
        \item $\tt{03431}$ has suffix $\tt{431}$, and $(\tt{4},\tt{3},\tt{1},\tt{3})$ is $\tfrac{16}{7}$-persistent with tipping point $5$,
        \item $\tt{03432}$ has suffix $\tt{432}$, and $(\tt{4},\tt{3},\tt{2},\tt{3})$ is $\tfrac{16}{7}$-persistent with tipping point $2$,
        \item $\tt{03433}$ has suffix $\tt{33}$, and $(\varepsilon,\tt{33},\varepsilon,\tt{3})$ is $\tfrac{16}{7}$-persistent with tipping point $14$,
        \item $\tt{03434}$ has suffix $434$, and $(\tt{4},\tt{3},\tt{4},\tt{3})$ is $\tfrac{16}{7}$-persistent with tipping point $1$,
        \item $\tt{0344}$ has suffix $344$, and $(\tt{3},\tt{44},\varepsilon,\tt{4})$ is $\tfrac{16}{7}$-persistent with tipping point $1$.
    \end{itemize}
    Thus, we see that the only possible complete return words to $\tt{0}$ with prefix $\tt{03}$ satisfying the given conditions are $\tt{030}$ and $\tt{03430}$, as reported in the proof of \autoref{Lemma:TernaryInductiveStep}.
\end{example}

\begin{example}\label{Example:PsiBacktracking}
    Consider the search for all words $w$ such that 
    \begin{enumerate}[label=(\roman*)]
        \item $\varphi^2(\psi_4(w))$ is rich, and 
        \item $w$ contains no bookended repetition $R$ such that $\psi_4(R)$ is $\tfrac{16}{7}$-persistent with tipping point at most $100$. 
    \end{enumerate}
    We illustrate the search for all such words in \autoref{Figure:PsiBacktracking}.  Each yellow leaf in the tree corresponds to a word $w$ that fails condition (i), i.e., $\varphi^2(\psi_4(w))$ is not rich.  In fact, each word corresponding to a yellow leaf in the tree has some word from the set $U=\{\tt{12}, \tt{13}, \tt{21}, \tt{31}\}$ as a suffix, and it is easy to verify by hand that for each word $u\in U$, the word $\varphi^2(\psi_4(u))$ is not rich.  We make special note here of the fact that the words $\varphi(\psi_4(\tt{12}))$ and $\varphi(\psi_4(\tt{21}))$ are rich!  This is one of the reasons that we assumed that $\varphi^2(\xx)$ is rich in the statement of \autoref{Lemma:TernaryInductiveStep} (instead of the weaker assumption that $\xx$ is rich), and why we used $\hat{\pi}$ in place of $\pi$ in the statement of \autoref{Theorem:TernaryStructureInduction}.
    Finally, each red leaf in the tree corresponds to a word that fails condition (ii).  In fact, each word corresponding to a red leaf in the tree has some word from the set
    \[
        F=\{\tt{00},\tt{0101},\tt{02},\tt{03},\tt{1010},\tt{11},\tt{20},\tt{22},\tt{233},\tt{30},\tt{323},\tt{332},\tt{333}\}
    \]
    as a suffix, and we explain below why each factor in $F$ fails condition (ii).%
    \begin{figure}
        \centering
        \begin{forest}
        [
        $\varepsilon$,for tree={grow=0,l=1.8cm}
            [
                $\tt{3}$
                [
                    $\tt{3}$
                    [
                        $\tt{3}$,draw,fill=red
                    ]
                    [
                        $\tt{2}$,draw,fill=red
                    ]
                    [
                        $\tt{1}$,draw,fill=yellow
                    ]
                    [
                        $\tt{0}$,draw,fill=red
                    ]
                ]
                [
                    $\tt{2}$
                    [
                        $\tt{3}$,draw,fill=red
                    ]
                    [
                        $\tt{2}$,draw,fill=red
                    ]
                    [
                        $\tt{1}$,draw,fill=yellow
                    ]
                    [
                        $\tt{0}$,draw,fill=red
                    ]
                ]
                [
                    $\tt{1}$,draw,fill=yellow
                ]
                [
                    $\tt{0}$,draw,fill=red
                ]
            ]
            [
                $\tt{2}$
                [
                    $\tt{3}$
                    [
                        $\tt{3}$,draw,fill=red
                    ]
                    [
                        $\tt{2}$
                        [
                            $\tt{3}$,draw,fill=red
                        ]
                        [
                            $\tt{2}$,draw,fill=red
                        ]
                        [
                            $\tt{1}$,draw,fill=yellow
                        ]
                        [
                            $\tt{0}$,draw,fill=red
                        ]
                    ]
                    [
                        $\tt{1}$,draw,fill=yellow
                    ]
                    [
                        $\tt{0}$,draw,fill=red
                    ]
                ]
                [
                    $\tt{2}$,draw,fill=red
                ]
                [
                    $\tt{1}$,draw,fill=yellow
                ]
                [
                    $\tt{0}$,draw,fill=red
                ]
            ]
            [
                $\tt{1}$
                [
                    $\tt{3}$,draw,fill=yellow
                ]
                [
                    $\tt{2}$,draw,fill=yellow
                ]
                [
                    $\tt{1}$,draw,fill=red
                ]
                [
                    $\tt{0}$
                    [
                        $\tt{3}$,draw,fill=red
                    ]
                    [
                        $\tt{2}$,draw,fill=red
                    ]
                    [
                        $\tt{1}$
                        [
                            $\tt{3}$,draw,fill=yellow
                        ]
                        [
                            $\tt{2}$,draw,fill=yellow
                        ]
                        [
                            $\tt{1}$,draw,fill=red
                        ]
                        [
                            $\tt{0}$,draw,fill=red
                        ]
                    ]
                    [
                        $\tt{0}$,draw,fill=red
                    ]
                ]
            ]
            [
                $\tt{0}$
                [
                    $\tt{3}$,draw,fill=red
                ]
                [
                    $\tt{2}$,draw,fill=red
                ] 
            	[
            		$\tt{1}$
                    [
                        $\tt{3}$,draw,fill=yellow            
                    ]
                    [
                        $\tt{2}$,draw,fill=yellow            
                    ]
                    [
                        $\tt{1}$,draw,fill=red
                    ]
                    [
                        $\tt{0}$
                        [
                            $\tt{3}$,draw,fill=red
                        ]
                        [
                            $\tt{2}$,draw,fill=red
                        ]
                        [
                            $\tt{1}$,draw,fill=red
                        ]
                        [
                            $\tt{0}$,draw,fill=red
                        ]
                    ]
            	]
            	[
            		$\tt{0}$,draw,fill=red
            	]
            ]
        ]
        \end{forest}
        \caption{An illustration of the backtracking search for words $w$ satisfying conditions (i) and (ii) of \autoref{Example:PsiBacktracking}.}
        \label{Figure:PsiBacktracking}
    \end{figure}
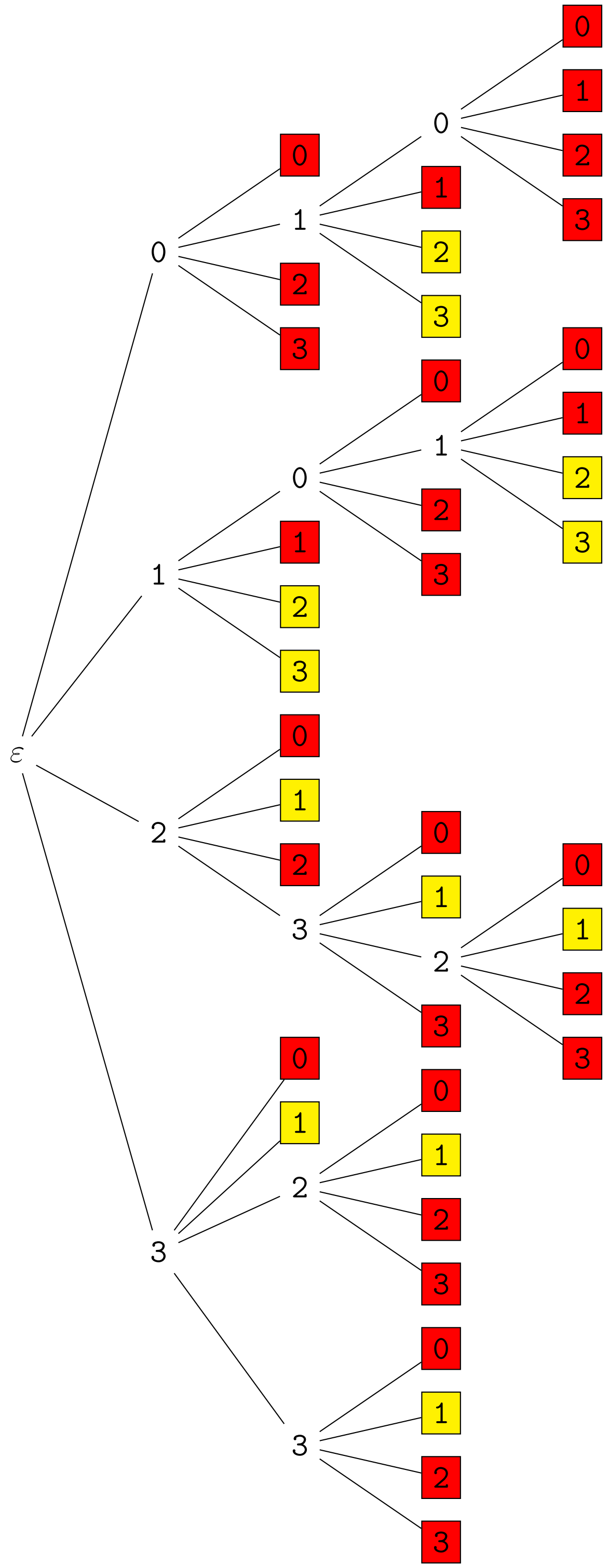%
    \begin{itemize}
        \item $\tt{00}$: $\psi_4(\varepsilon,\tt{00},\varepsilon,\tt{0})=(\varepsilon,\tt{01010},\varepsilon,\tt{01})$ is $\tfrac{16}{7}$-persistent with tipping point $4$,
        \item $\tt{0101}$: $\psi_4(\varepsilon,\tt{0101},\varepsilon,\tt{01})=(\varepsilon,\tt{010201020},\varepsilon,\tt{0102})$ is $\tfrac{16}{7}$-persistent with tipping point $7$,
        \item $\tt{02}$: $\psi_4(\varepsilon,\tt{0},\tt{2},\tt{0})=(\varepsilon,\tt{0101},\tt{3},\tt{01})$ is $\tfrac{16}{7}$-persistent with tipping point $4$,
        \item $\tt{03}$: $\psi_4(\varepsilon,\tt{0},\tt{3},\tt{0})=(\varepsilon,\tt{0101},\tt{4},\tt{01})$ is $\tfrac{16}{7}$-persistent with tipping point $4$,
        \item $\tt{1010}$: $\psi_4(\varepsilon,\tt{1010},\varepsilon,\tt{10})=(\varepsilon,\tt{020102010},\varepsilon,\tt{0201})$ is $\tfrac{16}{7}$-persistent with tipping point $7$,
        \item $\tt{11}$: $\psi_4(\varepsilon,\tt{11},\varepsilon,\tt{1})=(\varepsilon, \tt{02020},\varepsilon,\tt{02})$ is $\tfrac{16}{7}$-persistent with tipping point $5$,
        \item $\tt{20}$: $\psi_4(\tt{2},\tt{0},\varepsilon,\tt{0})=(\tt{3},\tt{1010},\varepsilon,\tt{10})$ is $\tfrac{16}{7}$-persistent with tipping point $4$, 
        \item $\tt{22}$: $\psi_4(\varepsilon,\tt{22},\varepsilon,\tt{2})=(\varepsilon,\tt{013101310},\varepsilon,\tt{0131})$ is $\tfrac{16}{7}$-persistent with tipping point $13$,
        \item $\tt{233}$: $\psi_4(\tt{2},\tt{33},\varepsilon,\tt{3})=(\tt{3},\tt{1014101410},\varepsilon,\tt{1014})$ is $\tfrac{16}{7}$-persistent with tipping point $1$,
        \item $\tt{30}$: $\psi_4(\tt{3},\tt{0},\varepsilon,\tt{0})=(\tt{4},\tt{1010},\varepsilon,\tt{10})$ is $\tfrac{16}{7}$-persistent with tipping point $4$,
        \item $\tt{323}$: $\psi_4(\tt{3},\tt{2},\tt{3},\tt{2})=(\tt{4},\tt{1013101},\tt{4},\tt{1013})$ is $\tfrac{16}{7}$-persistent with tipping point $3$,
        \item $\tt{332}$: $\psi_4(\varepsilon,\tt{33},\tt{2},\tt{3})=(\varepsilon,\tt{0141014101},\tt{3},\tt{0141})$ is $\tfrac{16}{7}$-persistent with tipping point $1$,
        \item $\tt{333}$: $\psi_4(\varepsilon,\tt{333},\varepsilon,\tt{3})=(\varepsilon,\tt{0141014101410},\varepsilon,\tt{0141})$ is $\tfrac{16}{7}$-persistent with tipping point $0$.
    \end{itemize}
    Thus, we see that the longest word over $\Sigma_4$ satisfying conditions (i) and (ii) has length $3$, as reported in \autoref{Table:PsiBacktracking} in the proof of \autoref{Lemma:TernaryInductiveStep}.
\end{example}

\end{document}